\documentclass[11pt]{amsart}

\newcommand{\FF}{\mathbb{F}}
\newcommand{\GG}{\mathbb{G}}
\newcommand{\NN}{\mathbb{N}}

\newcommand{\ZZ}{\mathbb{Z}}
\newcommand{\RR}{\mathbb{R}}
\newcommand{\St}{\mathrm{St}}

\newcommand{\cB}{\mathcal{B}}
\newcommand{\cC}{\mathcal{C}}

\newcommand{\cF}{\mathcal{F}}

\newcommand{\cN}{\mathcal{N}}

\newcommand{\cV}{\mathcal{V}}

\newcommand{\fg}{\mathfrak{g}}

\newcommand{\fu}{\mathfrak{u}}
\newcommand{\fv}{\mathfrak{v}}

\newcommand{\fF}{\mathfrak{F}}

\newcommand{\dact}{\boldsymbol{.}}
\newcommand{\lra}{\longrightarrow}

\DeclareMathOperator{\Ad}{Ad}

\DeclareMathOperator{\cx}{cx}
\DeclareMathOperator{\depth}{dp}
\DeclareMathOperator{\Dist}{Dist}
\DeclareMathOperator{\Ext}{Ext}
\DeclareMathOperator{\HH}{H}

\DeclareMathOperator{\Ht}{Ht}
\DeclareMathOperator{\Hom}{Hom}
\DeclareMathOperator{\id}{id}
\DeclareMathOperator{\im}{im}

\DeclareMathOperator{\Lie}{Lie}
\DeclareMathOperator{\modd}{mod}
\DeclareMathOperator{\Rad}{Rad}

\DeclareMathOperator{\SL}{SL}
\DeclareMathOperator{\Soc}{Soc}
\DeclareMathOperator{\supp}{supp}
\DeclareMathOperator{\Top}{Top}

\usepackage{pictex}
\usepackage{amscd}
\usepackage{latexsym}
\usepackage{amssymb}
\usepackage{eufrak}
\usepackage{verbatim}

\numberwithin{equation}{section}

\newtheorem{Theorem}[equation]{Theorem}
\newtheorem{Lemma}[equation]{Lemma}
\newtheorem{Proposition}[equation]{Proposition}
\newtheorem{Corollary}[equation]{Corollary}

\newtheorem*{thm*}{Theorem}
\theoremstyle{remark}
\newtheorem*{Remark}{Remark}
\newtheorem*{Remarks}{Remarks}

\newtheorem*{Example}{Example}

\numberwithin{equation}{section}

\begin{document}

\title[Varieties and Good Filtrations]{Support Varieties, AR-Components, and Good Filtrations}

\author[R. Farnsteiner \lowercase{and} G. R\"ohrle]{Rolf Farnsteiner \lowercase{and} Gerhard R\"ohrle}

\address[R. Farnsteiner]{Department of Mathematics, University of Kiel, Ludewig-Meyn-Str.~4, 24098 Kiel, Germany}
\email{rolf@math.uni-kiel.de}

\address[G. R\"{o}hrle]{Department of Mathematics, University of Bochum, 44780 Bochum, Germany}
\email{gerhard.roehrle@ruhr-uni-bochum.de}

\date{\today}

\makeatletter
\makeatother

\subjclass[2000]{Primary 16G70, Secondary 17B50}

\begin{abstract} Let $G$ be a reductive group, defined over the Galois field $\FF_p$ with $p$ being good for $G$. Using support varieties and covering techniques based on $G_rT$-modules,
we determine the position of simple modules and baby Verma modules within the stable Auslander-Reiten quiver $\Gamma_s(G_r)$ of the $r$-th Frobenius kernel of $G$. In particular, we
show that the almost split sequences terminating in these modules usually have an indecomposable middle term.  

Concerning support varieties, we introduce a reduction technique leading to isomorphisms 
\[\cV_{G_r}(Z_r(\lambda)) \cong \cV_{G_{r-d}}(Z_{r-d}(\mu))\] 
for baby Verma modules of certain highest weights $\lambda, \mu \in X(T)$, which are related by the notion of depth.  \end{abstract}

\maketitle

\setcounter{section}{-1}

\section{Introduction} \label{s:intro}
In the representation theory of finite-dimensional self-injective algebras, the stable Auslander--Reiten quiver has proven to be an important homological invariant, which has been studied for 
group algebras of finite groups, reduced enveloping algebras of restricted Lie algebras and distribution algebras of infinitesimal group schemes. In the classical context of Frobenius kernels of 
reductive groups, the relevant algebras are usually of wild representation type, rendering a classification of their indecomposable modules a hopeless task. This fact notwithstanding, one does 
have a fairly good understanding of the connected components of the corresponding stable Auslander--Reiten quiver. It is therefore of interest to relate this information to certain classes of 
indecomposable modules, such as simple modules, Weyl modules or Verma modules, and to identify the position of the latter within the AR-quiver. The main problem in this undertaking is the 
lack of a suitable presentation of the underlying algebras, as required by the techniques of abstract representation theory. On the other hand, the theory of rank varieties and support varieties 
has seen considerable progress over the last years, so that one can hope to exploit these tools in the aforementioned context. 

In continuation of work begun in \cite{FR}, we study in this article those Auslander-Reiten components of the algebras $\Dist(G_r)$ that contain simple modules or baby Verma modules. 
The underlying algebra $\Dist(G_r)$  consists of the distributions of the $r$-th Frobenius kernel of the smooth reductive group scheme $G$, defined over an algebraically closed field of positive
characteristic $p$. The case $r=1$, which pertains to restricted enveloping algebras of reductive Lie algebras, was settled in \cite{FR} by means of a detailed analysis of nilpotent orbits in rank 
varieties, leading to the consideration of groups of types $\SL(2)_1\!\times\!\SL(2)_1$, $\SL(3)_1$ and ${\rm SO}(5)_1$.  For $r>1$, rank varieties are less tractable and the structure of the
cohomology rings defining support varieties is more complicated. We address these problems by passage to $G_rT$-modules, as described below. With regard to Auslander-Reiten theory, our 
main results can roughly be summarized as follows:

\bigskip

\begin{thm*} Let $G$ be a smooth reductive group scheme, defined over $\FF_p$. Given a character $\lambda \in X(T)$ and $r \ge 1$, the following statements hold:

{\rm (1)} \ The simple $G_r$-module $L_r(\lambda)$ is either projective or it belongs to an $AR$-component of tree class $\tilde{A}_{12}$ or $A_\infty$. In the latter case, the middle term of the almost split sequence terminating in $L_r(\lambda)$ is indecomposable.

{\rm (2)} \  The baby Verma module $Z_r(\lambda)$ is either projective or it belongs to an $AR$-component of tree class $A_\infty$, with the middle term of the almost split sequence terminating in $Z_r(\lambda)$ being indecomposable. \end{thm*}

\bigskip
\noindent
The approach chosen in this article differs from that of \cite{FR} by the systematic use of Jantzen's category $\modd G_rT$ of $G_rT$-modules, which affords almost split sequences. Exploiting 
the natural ordering on the weights $X(T)$, we first study AR-components for $G_rT$-modules, and then use covering properties of the forgetful functor $\fF : \modd G_rT \lra \modd G_r$ to 
obtain information on the corresponding $G_r$-modules. The second major advantage of working in $\modd G_rT$ is its tractability with respect to coverings $\tilde{G} \lra G$, ultimately 
allowing us to bring Steinberg's tensor product theorem to bear. Aside from these technical aspects, the highest weight category $\modd G_rT$ possesses the much studied subcategory 
$\cF(\Delta)$ of $\widehat{Z}_r$-filtered modules. According to Ringel's seminal work \cite{Ri2}, the corresponding categories of $\Delta$-good modules over finite-dimensional  
quasi-hereditary algebras afford relative almost split sequences. It turns out that in our context these are merely the AR-sequences of the ambient Frobenius category $\modd G_rT$. 

Our paper is organized as follows. After a preliminary section, we determine in Section \ref{s:simples} the location of the simple $G_r$-modules. For finite groups and restricted enveloping algebras the corresponding problems were studied in \cite{Jo, Ka, KMU1, KMU2} and \cite{Fa2}, respectively. In fact, since the algebra $\Dist(G_r)$ is symmetric, Kawata's results \cite{Ka}, 
that were inspired by the modular representation theory of finite groups, can be brought to bear.  With the exception of simple modules belonging to blocks of tame representation type, the AR-components containing simple modules are of type $\ZZ[A_\infty]$, with the simple vertex having only one predecessor (see part (1) of the above Theorem). 

In Section \ref{s:AR-GT} we recall basic results from \cite{Fa5} concerning the AR-Theory of the Frobenius category $\modd G_rT$. As a first application, we exploit coverings $\tilde{G} \lra 
G$ to show that each AR-component of $\modd G_r$ contains at most one simple module. Using fundamental properties of the highest weight category $\modd G_rT$, we verify in Section 
\ref{s:AR-Ve} the second part of the above Theorem (cf.\ \cite[Thm.~B]{Un} for the corresponding result concerning trivial source modules over finite groups). 

Sections \ref{s:grvar} and \ref{s:depth} are concerned with filtrations and varieties of Verma modules. Motivated by Jantzen's work on blocks of $\Dist(G_r)$ \cite{Ja1}, we subdivide the defining highest weights of Verma modules according to their ``depth''. It turns out that questions concerning varieties and AR-components of Verma modules can often be reduced to the consideration of highest weights of depth $1$. Our method is based on an isomorphism
\[ \cV_{G_r}(M^{[d]}\!\otimes_k\!\St_d) \cong \cV_{G_{r-d}}(M),\]
relating the variety of a $G_{r-d}$-module $M$ to that of the tensor product $M^{[d]}\!\otimes_k\!\St_d$ of its $d$-th Frobenius twist $M^{[d]} \in \modd G_r$ with the $d$-th
Steinberg module $\St_d$. When applied to baby Verma modules, this yields explicit information on the varieties $\cV_{G_r}(Z_r(\lambda))$.

Our final Section \ref{s:appl} illustrates the utility of our results and techniques by providing applications concerning Verma modules of complexity at most $2$. For instance, the supports of
Verma modules of complexity $2$ are shown to be equidimensional.

\bigskip

\section{Preliminaries}\label{s:prelims}
Throughout, we will be working over an algebraically closed field $k$ of characteristic $p \ge 3$. Unless mentioned otherwise, all algebras and modules are assumed to be finite-dimensional.

Let $G$ be a connected, reductive, smooth group scheme over $k$ with maximal torus $T$ and character group $X(T)$. We fix a Borel subgroup $B = TU$ with unipotent radical $U$ and 
denote the Lie algebra of $U$ by $\fu = \Lie(U)$. Moreover, we write $\Psi = \Psi(T)$, $\Psi^+$ and $\Sigma$ for the set of roots of $G$ relative to $T$, the sets of positive and simple 
roots of $\Psi$ relative to $B$, respectively.

For $\lambda \in X(T)$  and $r \in \NN$ we put
\[ \Psi^r_\lambda := \{\alpha\in\Psi \mid \langle\lambda + \rho, \alpha^\vee \rangle \in p^r\mathbb Z\}, \]
where $\rho := \frac{1}{2} \sum_{\alpha \in \Psi^+}\alpha$ denotes the half-sum of the positive roots. We say that $\lambda$ is $p^r$-\emph{regular} (relative to $\Psi$) 
provided $\Psi^r_\lambda = \varnothing$. For $r=1$ we simply write $\Psi_\lambda = \Psi^1_\lambda$.
The set of $p^r$-\emph{restricted weights} $X_r(T)$ of $T$ is defined by
\[  X_r(T) := \{\lambda \in X(T) \mid 0 \le \langle \lambda , \alpha^\vee \rangle < p^r ; \alpha \in\Sigma\}.\]
For $G = \SL(2)$ we have $X(T) \cong \ZZ$ and thus identify weights with integers in this case.

Given $r \in \NN$, we let $F^r : G \lra G^{(r)}$ be the \emph{r-th Frobenius homomorphism\/} of $G$, whose comorphism is
\[ k[G]^{(r)} \lra k[G] \ ; \ x \mapsto x^{p^r}.\]
As usual, $k[G]^{(r)}$ denotes the $k$-algebra, whose underlying $\ZZ$-algebra structure is that given by $k[G]$, with $k$ acting via $\alpha\dact x := \alpha^{p^{-r}}x$ (cf.\ \cite[(I.9)]{Ja3}). The infinitesimal algebraic $k$-group $G_r := \ker F^r$  is called the \emph{r-th Frobenius kernel\/} of $G$.

For a $G_r$-module $M$ we define its \emph{cohomological support variety\/} $\cV_{G_r}(M)$ as the variety of the kernel of the canonical homomorphism
\[ \HH^{\rm ev}(G_r, k) \lra \Ext_{G_r}^{\rm ev}(M,M), \]
see \cite{SFB2}. We consider the $r$-th Frobenius kernel $\GG_{a(r)} := {\rm Spec}_k(k[T]/(T^{p^r}))$ of the additive group $\GG_a \cong {\rm Spec}_k(k[T])$ and recall that
\[ \Dist(\GG_{a(r)}) \cong  k[X_0, \ldots, X_{r-1}]/(X_0^p,\ldots,X_{r-1}^p),\]
with $X_i$ corresponding to the functional $x_i \in \Dist(\GG_{a(r)}) = k[\GG_{a(r)}]^\ast$ given by
\[ x_i(T^j + (T^{p^r})) = \delta_{p^ij}.\]
In particular, the algebra $\Dist(\GG_{a(r)}) = k[x_1,\ldots,x_r]$ is generated by the $x_i$. Owing to \cite[(6.8)]{SFB2}, the support variety of $M$ is homeomorphic to the \emph{rank 
variety}
\[ V_r(G)_M := \{ \varphi : \GG_{a(r)} \lra  G \mid M|_{A_r} \ {\rm is \ not \ projective}\}.\]
Here $M|_{A_r}$ denotes the pull-back of the module structure of $M$ along the restriction to $A_r := k[x_{r-1}] $ of the homomorphism $\Dist(\GG_{a(r)}) \lra \Dist(G_r)$ corresponding 
to $\varphi$. In \cite{SFB2} the authors consider the support scheme $V_r(G)_M$, whose $k$-rational points are our variety $V_r(G)_M$. 

By general theory, the category of $\modd G_r$ of $G_r$-modules is equivalent to the module category $\modd \Dist(G_r)$ of the Hopf algebra $\Dist(G_r)$ of distributions of $G_r$, see 
\cite[I, \S8]{Ja3}. According to \cite[(II, \S7, 4.2)]{DG}, the Hopf algebra $\Dist(G_1)$ is isomorphic to the {\it restricted enveloping algebra} $U_0(\fg)$ of the restricted Lie algebra 
$\fg = \Lie(G)$ of $G$. This isomorphism induces an isomorphism between $\cV_{G_1}(M)$ and the corresponding variety $\cV_{\fg}(M)$ of the $\fg$-module associated to $M$.
We shall henceforth identify these two varieties without further notice.

According to \cite[(5.4)]{SFB2}, a closed embedding $H \hookrightarrow H'$ of infinitesimal group schemes gives rise to an embedding $\cV_H(k) \hookrightarrow \cV_{H'}(k)$ that maps 
$\cV_H(k)$ homeomorphically onto its image. We use this identification and write
\[ \cV_H(M) = \cV_{H'}(M) \cap \cV_H(k)\]
for an $H'$-module $M$. For future reference we record the following well-known fact concerning varieties of relatively projective modules, see for instance \cite[(2.3.1)]{NPV}.

\bigskip

\begin{Proposition} \label{Pr1} Let $H \subseteq H'$ be infinitesimal group schemes. 

{\rm (1)} \ If $M$ is a $(\Dist(H')\!:\!\Dist(H))$-projective $H$-module, then $\cV_{H'}(M) = \cV_{H}(M)$.

{\rm (2)} \ If $N$ is an $H$-module, then $\cV_{H'}(\Dist(H')\!\otimes_{\Dist(H)}\!N)\subseteq \cV_H(N)$. 

{\rm (3)} \ A $\Dist(H)$-module $M$ is projective if and only if $\cV_H(M) = \{0\}$. \hfill $\square$ \end{Proposition}

\bigskip
\noindent
We let 
\[X(T)_+ := \{ \lambda \in X(T) \mid \langle \lambda, \alpha^\vee \rangle \ge 0 \ \ \ \forall \ \alpha \in \Sigma\}\] 
be the set of \emph{dominant weights} and recall that for every $\lambda \in X(T)_+$ there exists a simple $G$-module $L(\lambda)$ of highest weight $\lambda$. Moreover, the 
$L(\lambda)$ form a complete set of representatives for the isoclasses of the simple $G$-modules. Similarly, the simple $G_r$-modules are of the form $L_r(\lambda)$, with $\lambda$ 
belonging to a set of representatives for $X(T)/p^rX(T)$ (cf.\ \cite[(II.3.10)]{Ja3}). Finally, the simple $G_rT$-modules are given by $\widehat{L}_r(\lambda)$, with $\lambda \in X(T)$. 
Note that $\widehat{L}_r(\lambda)|_{G_r} \cong L_r(\lambda)$ (see \cite[(II.9.6)]{Ja3}).

Given $\lambda \in X(T)$ and $r \in \NN$, we let
\[  Z_r(\lambda) := \Dist(G_r)\!\otimes_{\Dist(B_r)}\!k_\lambda \]
be the (baby) \emph{Verma module\/} of $G_r$ with highest weight $\lambda$.\footnote{Our notation differs from that in Jantzen's book \cite{Ja3}, whose notational conventions we 
follow fairly closely otherwise. In \cite{Ja3} our group $B$ is denoted $B^+$.} For a Levi subgroup $L \supseteq T$ of $G$ we define the  (baby) \emph{Verma module\/} of $L_r$ with 
highest weight $\lambda$ by
\[  Z_r^L(\lambda) := \Dist(L_r)\!\otimes_{\Dist((B\cap L)_r)}\!k_\lambda. \]
We record the following basic result concerning varieties (cf.\ \cite[(2.3.1),(4.2.1)]{NPV}):

\bigskip

\begin{Proposition} \label{Pr2} Let $\lambda \in X(T)$. Then the following statements hold:

{\rm (1)} \ $\cV_{G_r}(Z_r(\lambda)) \subseteq \cV_{U_r}(k)$.

{\rm (2)} \ $\cV_{L_r}(Z_r^L(\lambda)) \subseteq \cV_{G_r}(Z_r(\lambda))$. \hfill $\square$ \end{Proposition}

\bigskip
\noindent
We denote by $W$ and $W_p$ the Weyl group and the affine Weyl group associated to the reductive group scheme $G$, respectively. The ``dot'' action of $w \in W_p$ on $X(T)$ is defined as 
follows:
\[ w \dact \lambda := w(\lambda + \rho) -\rho.\]
Let $H$ be an infinitesimal group scheme,  $M$ be an $H$-module. By definition, the \emph{complexity}  $\cx_H(M)$ of $M$ is the polynomial rate of growth of a minimal projective 
resolution $(P_i)_{i\ge 0}$ of $M$, i.e., 
\[ \cx_H(M) := \min\{ n \in \NN_0 \cup \{\infty\}  \mid \exists \ c > 0 \ \dim_k P_i \le c i^{n-1} \text{ for all } i \ge 1\}. \]
The reader is referred to \cite[(5.3)]{Be2} for basic properties of this notion.

Let $h = h(G)$ denote the Coxeter number of $G$, that is, the maximum of the Coxeter numbers of the simple components of the derived group $(G,G)$. 

Recall that a prime $p$ is said to be \emph{good\/} for $G$, provided it does not divide any of the coefficients occurring when expressing any root of $G$ as a linear combination of simple 
roots. Owing to \cite[(2.7)]{Ja2}, $\Psi^r_\lambda$ is a subsystem of $\Psi$ whenever $p$ is good for $G$. 

\bigskip

\section{Auslander-Reiten Components of Simple $G_r$-modules} \label{s:simples}
Given a self-injective algebra $\Lambda$, we denote by $\Gamma_s(\Lambda)$ the {\it stable Auslander--Reiten quiver} of $\Lambda$. By definition, the directed graph 
$\Gamma_s(\Lambda)$ has as vertices the non-projective indecomposable $\Lambda$-modules and its arrows are defined via the so-called {\it irreducible morphisms}. We refer the interested 
reader to \cite[Chap.\ VII]{ARS} for further details. The AR-quiver is fitted with an automorphism $\tau_\Lambda$, the so-called {\it Auslander--Reiten translation}. Since $\Lambda$ is 
self-injective, $\tau_\Lambda$ coincides with the composite $\Omega^2_\Lambda \circ \nu_\Lambda$ of the square of the Heller translate $\Omega_\Lambda$ and the Nakayama functor 
$\nu_\Lambda$, cf.\ \cite[(IV.3.7)]{ARS}.

The connected components of $\Gamma_s(\Lambda)$ are connected stable translation quivers. By work of Riedtmann \cite[Struktursatz]{Ri}, the structure of such a quiver $\Theta$ is 
determined by a directed tree $T_\Theta$ and an {\it admissible group} $\Pi \subseteq {\rm Aut}_k(\ZZ[T_\Theta])$, giving rise to an isomorphism
\[ \Theta \cong \ZZ[T_\Theta]/\Pi\]
of stable translation quivers. The undirected tree $\bar{T}_\Theta$ of $T_\Theta$ is uniquely determined by $\Theta$ and is called the {\it tree class} of $\Theta$. We refer the reader to 
\cite[(4.15.6)]{Be} for further details. For group algebras of finite groups, the possible tree classes and admissible groups were first determined by Webb \cite{We}, with refinements provided 
in \cite{Re,Ok,Bes1,Bes2,ES}.

{\it Throughout, $G$ is assumed to be a connected smooth reductive algebraic group scheme}. By work of Larson and Sweedler \cite{LS}, every finite-dimensional Hopf algebra is a Frobenius 
algebra. Consequently, the algebra $\Dist(G_r)$ is self-injective, so that the classes of projective and injective $G_r$-modules coincide.

Following Ringel \cite{Ri1}, an indecomposable $G_r$-module $M$ is called \emph{quasi-simple} provided its isomorphism class $[M]$ lies at the end of a component of tree class $A_\infty$ 
of the stable AR-quiver $\Gamma_s(G_r)$ associated to $\Dist(G_r)$. Given $\lambda \in X(T)$, let $P_r(\lambda)$ denote the projective cover of the non-projective simple $G_r$-module 
$L_r(\lambda)$ and let $\Ht_r(\lambda) = \Rad(P_r(\lambda))/\Soc(P_r(\lambda))$ be its \emph{heart}. We write $\Top(M) := M/\Rad(M)$ for the {\it top} of a $G_r$-module $M$.

For future reference we record the following basic property of the algebras of distributions (see \cite{Sc,Hu1} for related work):

\bigskip

\begin{Lemma} \label{AR-S1} The Hopf algebra $\Dist(G_r)$ is symmetric. \end{Lemma}

\begin{proof} According to \cite[(I.9.7)]{Ja3}, the smooth group scheme $G$ acts on the space of right integrals of $\Dist(G_r)$ via the character $g \mapsto \det \Ad(g)^{p^r-1}$. Since 
$G = (G,G)Z(G)$, cf.\ \cite[(27.5)]{Hu}, this character is trivial, proving that the modular function of the Hopf algebra $\Dist(G_r)$ coincides with its counit. Consequently, $\Dist(G_r)$ is 
symmetric, see \cite[(1.5)]{FMS} for more details. \end{proof}

\bigskip
\noindent
For future reference, we reformulate parts of \cite[(7.1)]{Fa3}.

\bigskip

\begin{Lemma} \label{AR-S2} Let $\lambda \in X(T)$ be a weight. Then the following statements hold:

{\rm (1)} \ $\dim \cV_{G_r}(L_r(\lambda))\ne 1$.

{\rm (2)} \ If $\dim \cV_{G_r}(L_r(\lambda)) = 2$, then the block $\cB \subseteq \Dist(G_r)$ containing $L_r(\lambda)$ is Morita equivalent to a block of $\SL(2)_r$. \end{Lemma}

\begin{proof} (1) This was shown in the second paragraph of \cite[p.80]{Fa3}.

(2) This follows directly from \cite[p.80]{Fa3}. \end{proof}

\bigskip

\begin{Theorem} \label{AR-S3} Let $\lambda \in X(T)$ be a weight such that $L_r(\lambda)$ is not projective.

{\rm (1)} \ If $\dim \cV_{G_r}(L_r(\lambda)) \ne 2$, then $L_r(\lambda)$ is quasi-simple.

{\rm (2)} \ If $\dim \cV_{G_r}(L_r(\lambda)) = 2$, then either $L_r(\lambda)$ is quasi-simple or $L_r(\lambda)$ belongs to a component of type $\ZZ [{\tilde A}_{12}]$. \end{Theorem}

\begin{proof} If $\dim \cV_{G_r}(L_r(\lambda)) \ge 3$, then, thanks to \cite[(2.2)]{Fa4}, the component $\Theta$ containing $[L_r(\lambda)]$ is isomorphic to $\ZZ [A_\infty]$.

Assuming $\Theta \cong \ZZ[A_\infty]$, we show that the isoclass $[L_r(\lambda)]$ is located at an end of $\Theta$. Suppose that $L_r(\lambda)$ is not quasi-simple. Since $\Dist(G_r)$ is 
symmetric, \cite[(1.5)]{Ka} provides simple $\Dist(G_r)$-modules $L_r(\mu_i) \not \cong L_r(\lambda)$ for $1 \le i \le n$ such that each projective cover $P_r(\mu_i)$ is uniserial, of length 
$\ell(P_r(\mu_i))=n+2$, with $\Top(\Rad(P_r(\mu_i))) \cong L_r(\mu_{i-1})$, where $\mu_0 = \lambda$. By the same token, $L_r(\lambda)$ has multiplicity $1$ in $P_r(\mu_i)$ for 
each $i$. 

Assuming $n \ge 2$, we thus have 
\[ \dim_k \Ext^1_{G_r}(L_r(\mu_1), L_r(\mu_2)) = 0 \quad \text{ and} \quad \dim_k \Ext^1_{G_r}(L_r(\mu_2), L_r(\mu_1)) = 1, \]
which contradicts \cite[(II.9.19(2))]{Ja3}. Consequently, $n = 1$, so that $\ell(P_r(\mu_1)) = 3$. Owing to \cite[(II.11.4)]{Ja3}, the module $P_r(\mu_1)$ has a filtration by baby Verma modules. Since $\ell(P_r(\mu_1)) = 3$, one $Z_r$-filtration factor, $Z_r(\gamma)$ say, is simple or projective. In view of \cite[(II.11.8)]{Ja3}, the factor $Z_r(\gamma)$ is simple and projective, so that the block containing $L_r(\lambda)$ is simple and $\dim \cV_{G_r}(L_r(\lambda)) = 0$, a contradiction.

By Lemma \ref{AR-S2}, the variety $\cV_{G_r}(L_r(\lambda))$ is not one-dimensional and so the proof of part (1) is complete.

Suppose that $\dim \cV_{G_r}(L_r(\lambda)) = 2$ and that $\Theta \not \cong \ZZ[A_\infty]$. By Theorem 4.1 in \cite{Fa4} and \cite{Fa7}, the component containing $[L_r(\lambda)]$ is 
of type $\ZZ [{\tilde A}_{12}]$. \end{proof}

\bigskip

\begin{Remark} For $r=1$, a simple $G_r$-module of complexity $2$ belongs to a component of type $\ZZ[\tilde{A}_{12}]$, cf.\ \cite[(5.2)]{Fa3}.  By contrast, simple modules
of complexity $2$ may belong to components of type $\ZZ[A_\infty]$, whenever $r \ge 2$ (cf.\ \cite{Fa7}). \end{Remark}

\bigskip
\noindent
We record an immediate consequence concerning the structure of the hearts of the principal indecomposable modules:

\bigskip

\begin{Corollary} \label{AR-S4} Let $\lambda \in X(T)$ be a weight such that $L_r(\lambda)$ is not projective.

{\rm (1)} \ If $\dim \cV_{G_r}(L_r(\lambda)) \ne 2$, then $\Ht_r(\lambda)$ is indecomposable.

{\rm (2)} \ If $\dim \cV_{G_r}(L_r(\lambda)) = 2$, then either $\Ht_r(\lambda)$ is indecomposable or $\Ht_r(\lambda) \cong L_r(\mu) \oplus L_r(\mu)$, where $L_r(\mu) \not \cong 
L_r(\lambda)$. \end{Corollary}

\begin{proof} Owing to Lemma \ref{AR-S1}, the algebra $\Dist(G_r)$ is symmetric, so that $\Soc(P_r(\lambda)) \cong L_r(\lambda)$. Hence
\[ (0) \lra \Rad(P_r(\lambda)) \lra P_r(\lambda) \oplus \Ht_r(\lambda) \lra P_r(\lambda)/L_r(\lambda) \lra (0)\]
is the standard almost split sequence originating in $\Rad(P_r(\lambda))$, see \cite[(V.5.5)]{ARS}. Since $\Omega_{G_r}$ is an auto-equivalence of the stable module category of 
$\Dist(G_r)$, \cite[(IV.3.5)]{ARS}, it induces an automorphism of $\Gamma_s(G_r)$, \cite[(X.1.9)]{ARS}. Consequently, $L_r(\lambda)$ and $P_r(\lambda)/L_r(\lambda) \cong 
\Omega^{-1}_{G_r}(L_r(\lambda))$ have the same number of non-projective predecessors. If the component $\Theta$ containing $L_r(\lambda)$ is isomorphic to $\ZZ[A_\infty]$, then the 
result thus follows from Theorem \ref{AR-S3}. Alternatively, $\dim \cV_{G_r}(L_r(\lambda))=2$ and $L_r(\lambda)$ belongs to a component of type $\ZZ[\tilde{A}_{12}]$. Thus, 
\cite[(IV.3.8.3)]{Er} implies that $\Ht_r(\lambda) \cong L_r(\mu)\oplus L_r(\mu)$. The assumption $L_r(\lambda) \cong L_r(\mu)$ yields $\dim_k 
\Ext^1_{G_r}(L_r(\lambda),L_r(\lambda)) = 2$. However, by Lemma \ref{AR-S2}, the block containing $L_r(\lambda)$ is Morita equivalent to a block of $\Dist(\SL(2)_r)$, and 
the simple modules of the latter algebra are known to afford no non-trivial self-extensions (cf.\ \cite[Thm.]{Pf}). \end{proof}

\bigskip

\begin{Remark} The exceptional case of decomposable hearts corresponds to those blocks $\cB_r(\lambda) \subseteq \Dist(G_r)$ that have tame representation type. The reader is referred to 
\cite[(I.4)]{Er} for background information on representation type. 

Suppose that $\Ht_r(\lambda)$ is decomposable. In light of \cite[(IV.3.8.3)]{Er} and (2) of Corollary \ref{AR-S4}, the algebra $\cB_r(\lambda)/\Soc(\cB_r(\lambda))$ is special biserial, and 
hence in particular tame or representation-finite (cf.\ \cite[(II.3.1)]{Er}). Since $\cB_r(\lambda)$ and $\cB_r(\lambda)/\Soc(\cB_r(\lambda))$ have the same non-projective indecomposable 
modules, we conclude that $\cB_r(\lambda)$ enjoys the same property. As $\dim \cV_{G_r}(L_r(\lambda)) = 2$, it follows that $\cB_r(\lambda)$ is tame. Conversely, if $L_r(\lambda)$ 
belongs to a tame block $\cB_r(\lambda)$ of $\Dist(G_r)$, then $\cB_r(\lambda)$ is Morita equivalent to a tame block of $\Dist(\SL(2)_1)$, \cite[(7.1)]{Fa3}. Consequently, the component 
containing $L_r(\lambda)$ is isomorphic to $\ZZ[\tilde{A}_{12}]$ and $\Ht_r(\lambda)$ is decomposable. \end{Remark}
 
\bigskip
\noindent
The final result of this section shows that, provided $G$ is defined over $\FF_p$, every AR-component contains at most one simple module. The technical hypothesis on the derived subgroup 
$(G,G)$ of $G$ will later be removed by means of covering techniques, see Theorem \ref{AR-GT6}. Since $G$ is defined over $\FF_p$, we have the Frobenius endomorphism $F : G \lra G$, 
and for every $G$-module $M$, there is the Frobenius twist $M^{[1]}$, obtained by composing the representation afforded by $M$ with $F$ (see \cite[(II.3.16)]{Ja3}). For future reference,
we record the following:

\bigskip

\begin{Lemma} \label{AR-S5} Let $\Theta \subseteq \Gamma_s(G_r)$ be a component which is isomorphic to $\ZZ[\tilde{A}_{12}]$. Then $\Theta$ contains exactly one simple module.
\end{Lemma}

\begin{proof} According to \cite[Thm.A]{We} the component $\Omega_{G_r}^{-1}(\Theta) \cong \Theta$ is attached to a principal indecomposable module. In view of \cite[(V.5.5)]{ARS} 
and Lemma \ref{AR-S1}, there thus exists a simple $G_r$-module $S$ such that $[\Omega^{-1}_{G_r}(S)] \in \Omega^{-1}_{G_r}(\Theta)$. As a result, $[S] \in \Theta$.

Corollary \ref{AR-S4} in conjunction with \cite[(IV.3.8.3)]{Er} implies that the block of $\Dist(G_r)$ containing $S$ possesses two simple modules $S$ and $T$, with
\[ (0) \lra \Rad(P(S)) \lra P(S)\oplus T\oplus T \lra P(S)/S \lra (0)\]
being the almost split sequence involving the projective cover $P(S)$ of $S$. In particular, $T$ belongs to the component $\Omega^{-1}_{G_r}(\Theta)$. By virtue of  \cite[(IV.3.8.3)]{Er}, 
we have $\Theta \neq \Omega_{G_r}(\Theta)$, so that $S$ is the only simple module belonging to $\Theta$. \end{proof}

\bigskip

\begin{Theorem} \label{AR-S6}  Suppose that $G$ is defined over $\FF_p$ and that the derived subgroup $(G,G)$ of $G$ is simply connected. If $S$ is a simple, non-projective 
$G_r$-module, then $[S]$ is the only simple vertex in its stable AR-component. \end{Theorem}

\begin{proof} Let $\Theta \subseteq \Gamma_s(G_r)$ be the stable AR-component containing the vertex $[S]$. According to Theorem \ref{AR-S3}, the component $\Theta$ is of type $\ZZ[A_\infty]$ or $\ZZ[\tilde{A}_{12}]$. In the latter case, the assertion follows from Lemma \ref{AR-S5}. 

Assuming $\Theta \cong \ZZ[A_\infty]$, we proceed by induction on $r$. If $r= 1$, then \cite[(4.1)]{Fa2} proves our claim. Now assume that $r\ge 2$ and let $T \not \cong S$ be another 
simple $G_r$-module which belongs to $\Theta$. According to Theorem \ref{AR-S3}, $S$ and $T$ lie at the end of $\Theta$. Thanks to Lemma \ref{AR-S1},  the Auslander--Reiten translate of 
$\Gamma_s(G_r)$ coincides with $\Omega^2_{G_r}$. Without loss of generality, there thus exists $n \in \NN$ with
\begin{equation} \label{e:1} 
T \cong \Omega^{2n}_{G_r}(S).
\end{equation}
As $(G,G)$ is simply connected, each element $\lambda \in X(T)$ can be written as $\lambda = \lambda_0 + p^r\lambda_1$, with $\lambda_0 \in X_r(T)$ and $\lambda_1 \in X(T)$ (cf.\ 
\cite[(II.9.14)]{Ja3}). Consequently, $X_r(T)$ contains a complete set of representatives for $X(T)/p^rX(T)$, and  there exist $\lambda, \mu \in X_r(T)$ such that $S = L_r(\lambda), T = L_r(\mu)$, with both modules being restrictions of simple $G$-modules, see \cite[(II.3.15)]{Ja3}. By the same token, there are decompositions $\mu = \mu_0+p\mu_1$ and $\lambda = \lambda_0+p\lambda_1$, where $\mu_0, \lambda_0 \in X_1(T)$ and $\mu_1,\lambda_1 \in X_{r-1}(T)$. In view of \cite[(II.3.16)]{Ja3}, we thus have
\[ L_r(\mu) \cong L_1(\mu_0)\!\otimes_k\!L_{r-1}(\mu_1)^{[1]} \ \ \text{and} \ \  L_r(\lambda) \cong L_1(\lambda_0)\!\otimes_k\!L_{r-1}(\lambda_1)^{[1]},\]
so that
\[ L_r(\mu)|_{G_1} = m_\mu L_1(\mu_0) \ \ \text{and} \ \  L_r(\lambda)|_{G_1} = m_\lambda L_1(\lambda_0).\]
Identity \eqref{e:1} in conjunction with standard properties of the Heller operator now yields:
\begin{equation} \label{e:2} 
m_\mu L_1(\mu_0) \cong  m_\lambda \Omega_{G_1}^{2n}(L_1(\lambda_0)) \oplus ({\rm proj.}).
\end{equation}
In particular,  the indecomposable $G_1$-module $\Omega_{G_1}^{2n}(L_1(\lambda_0))$ (cf.\ \cite{He}) is a direct summand of $L_1(\mu_0)$. Consequently, \cite[(2.5.4)]{Be} gives
\[ \dim_k \Ext^{2n}_{G_1}(L_1(\lambda_0),L_1(\lambda_0)) \le \dim_k \Hom_{G_1}(L_1(\mu_0),L_1(\lambda_0)) \le 1.\]
Owing to \cite[(2.1)]{Fa1}, this implies $\cx_{G_1}(L_1(\lambda_0)) \le 1$, which, in view of $G$ being reductive, entails the projectivity of $L_1(\lambda_0)$, see \cite[(3.1)]{Fa2}. By \eqref{e:2}, the simple $G_1$-module $L_1(\mu_0)$ also has this property. 

Let $e_S, e_T \in \Dist(G_1)$ be the primitive central idempotents defining the simple blocks containing $S|_{G_1}$ and $T|_{G_1}$, respectively. Since these are fixed by the adjoint 
representation of the connected group $G$, it follows that $e_S,e_T$ are also central idempotents of $\Dist(G_r) \supseteq \Dist(G_1)$ (see \cite[(II.10.3)]{Ja3} for more details). As $S$ 
and $T$ belong to the same block of $\Dist(G_r)$, the idempotents $e_S$ and $e_T$ act via the identity on $S$ and $T$. Consequently, $e_S = e_T$, so that the simple projective 
$G_1$-modules $L_1(\lambda_0)$ and $L_1(\mu_0)$ belong to the same block, whence $\Hom_{G_1}(L_1(\lambda_0),L_1(\mu_0)) = k$. Let $\gamma \in X(G)$ be the character 
giving the action of $G$ on $\Hom_{G_1}(L_1(\lambda_0),L_1(\mu_0))$. There results an isomorphism
\[ L_1(\mu_0) \cong L_1(\lambda_0)\!\otimes_k\!k_\gamma\]
of $G_r$-modules. As a result, we have
\[ T \cong L_1(\lambda_0)\!\otimes_k\!L_{r-1}(\mu_1)^{[1]}\!\otimes_k\!k_{-\gamma}.\]
Since $G_1 \subseteq \ker \gamma$, it follows that $\gamma \in pX(G)$ (cf.\ \cite[(II.1.18,II.3.7)]{Ja3}), so that the last two tensor factors are the Frobenius twist of 
$L_{r-1}(\mu_1)\!\otimes_k\!k_{\omega}$, for some $\omega \in X(G)$.

Owing to \cite[(II.10.5)]{Ja3}, the functor
\[ M \mapsto L_1(\lambda_0)\!\otimes_k\!M^{[1]}\]
induces an equivalence between $\modd G_{r-1}$ and the sum of those blocks of $\modd G_r$ given by $e_S$. As a result, the modules $L_{r-1}(\mu_1)\!\otimes_k\!k_\omega $ and 
$L_{r-1}(\lambda_1)$ belong to the same component of $\Gamma_s(G_{r-1})$. The inductive hypothesis now provides an isomorphism $L_{r-1}(\mu_1)\!\otimes_k\!k_\omega \cong 
L_{r-1}(\lambda_1)$ of $G_{r-1}$-modules, so that 
\[ T \cong L_1(\lambda_0)\!\otimes_k\!(L_{r-1}(\mu_1)\!\otimes_k\!k_\omega)^{[1]} \cong L_1(\lambda_0)\!\otimes_k\!L_{r-1}(\lambda_1)^{[1]} \cong S,\]
a contradiction. \end{proof}

\bigskip

\section{Auslander-Reiten Theory for $\modd G_rT$} \label{s:AR-GT}
Let $G$ be a connected, reductive, smooth group scheme with maximal torus $T$ and character group $X(T)$. Given $\lambda \in X(T)$, we denote by $\Theta_r(\lambda)$ the connected 
component of the stable Auslander--Reiten quiver $\Gamma_s(G_r)$ of $\Dist(G_r)$ containing the isoclass $[Z_r(\lambda)]$. Our main tool in this section is covering theory, appearing in the 
guise of Jantzen's category $\modd G_rT$ of $G_rT$-modules (cf.\ \cite[\S2]{Ja1}). It turns out that certain properties of $\Theta_r(\lambda)$ can be effectively studied by first investigating 
the corresponding properties for components of the stable Auslander--Reiten quiver of the Frobenius category $\modd G_rT$. In view of the conjugacy of maximal tori, our considerations do 
not depend on the choice of $T$. The category $\modd G_rT$ is a highest weight category in the sense of \cite{CPS}, whose standard objects $\Delta(\lambda)$ are the baby Verma modules 
$\widehat{Z}_r(\lambda)$, where $\lambda \in X(T)$. Recall that the $G_rT$-modules $\widehat{Z}_r(\lambda)$ and $\widehat{Z}'_r(\lambda)$ are given by
\[\widehat{Z}_r(\lambda) := \Dist(G_r)\!\otimes_{\Dist(B_r)}\! k_\lambda \ \ \text{and} \ \ \widehat{Z}'_r(\lambda) := \Hom_{\Dist(B^-_r)}(\Dist(G_r),k_\lambda),\]
where $B^-$ is the Borel subalgebra opposite to our standard Borel subalgebra $B \subseteq G$. The $T$-action is induced by the adjoint action of $T$ on $\Dist(G_r)$ (cf.\ \cite[\S 
II.9]{Ja3}). 

We shall consider the full subcategory $\cF(\Delta)$ of  $\modd G_rT$ of $\widehat{Z}_r$-filtered $G_rT$-modules. Given $M \in \cF(\Delta)$, its filtration multiplicities
\[ [M\!:\!\widehat{Z}_r(\lambda)]\]
are well-defined (cf.\ \cite[(II.11.2)]{Ja3}). We define the {\it $\Delta$-support} of $M \in \cF(\Delta)$ via
\[ \supp_\Delta(M) := \{ \lambda \in X(T) \mid [M\!:\!\widehat{Z}_r(\lambda)] \ne 0\}.\]
The objects of $\cF(\Delta)$ are usually referred to as {\it $\Delta$-good} modules. 

\bigskip
\noindent
Given $\lambda \in X(T)$, we write $X(T)_{\ge \lambda} :=\{ \mu \in X(T) \mid \mu \ge \lambda\}$ as well as $X(T)_{> \lambda} :=\{ \mu \in X(T) \mid \mu > \lambda\}$.

Recall that for every $\lambda \in X(T)$ there exists a unique simple object $\widehat{L}_r(\lambda)$ with projective cover $\widehat{P}_r(\lambda)$, and that all simple and projective 
indecomposable objects are of this form (see \cite[(II.9.6), (II.11.5.(3))]{Ja3}). In view of \cite[(II.11.4)]{Ja3}, the module $\widehat{P}_r(\lambda)$ belongs to $\cF(\Delta)$, with its 
filtration multiplicities obeying the BGG reciprocity formula
\[ [\widehat{P}_r(\lambda)\!:\!\widehat{Z}_r(\mu)] = [\widehat{Z}_r(\mu)\!:\!\widehat{L}_r(\lambda)].\]
In particular, the highest weights of the $\widehat{Z}_r$-filtration factors of $\widehat{P}_r(\lambda)$ enjoy the following property:

\bigskip

\begin{Lemma} \label{AR-GT1} We have
\[ \supp_\Delta(\widehat{P}_r(\lambda)) \subseteq X(T)_{\ge \lambda}\]
for every $\lambda \in X(T)$. \hfill $\square$ \end{Lemma}

\bigskip
\noindent
In the sequel, we shall study $\modd G_rT$ as well as $\modd (G_r\!\rtimes\!T)$. The latter category coincides with the category $\modd_{X(T)}\Dist(G_r)$ of $X(T)$-graded $G_r$-modules 
and degree zero homomorphisms. We identify $X(T)$ with the subgroup of those characters $\lambda \in X(G_r\!\rtimes\!T)$ that are trivial on $G_r$. Thus, every $\lambda \in X(T)$ 
defines a one-dimensional $(G_r\!\rtimes\!T)$-module $k_\lambda$. It now follows directly from the definition that the shift functor $M \mapsto M\langle \lambda \rangle$ of 
$\modd_{X(T)}\Dist(G_r)$ corresponds to the auto-equivalence $M \mapsto M\!\otimes_k\!k_\lambda$.

Since $G_rT \cong (G_r\rtimes T)/T_r$, the shifts sending $\modd G_rT$ onto itself are given by those $\lambda \in X(T)$ that vanish on $T_r$. In view of \cite[(II.3.7)]{Ja3} there exists an 
exact sequence
\[ (0) \lra p^rX(T) \stackrel{\rm can.}{\lra} X(T) \stackrel{\rm res.}{\lra} X(T_r) \lra (0),\]
so that the relevant shifts are those belonging to the subgroup $p^rX(T)$ of $X(T)$.

The block decomposition of the algebraic group $G_r\!\rtimes\!T$ given in \cite[(II.7.1)]{Ja3} yields a direct sum decomposition
\[ \modd (G_r\!\rtimes\!T) = \bigoplus_{b \in \cB(G_r\!\rtimes\!T)} (\modd (G_r\!\rtimes\!T))_b,\]
whose constituents are referred to as {\it blocks}. The proof of \cite[(2.1)]{Fa5} for $r=1$ may be adopted verbatim to obtain:

\bigskip

\begin{Lemma} \label{AR-GT2} The following statements hold:

{\rm (1)} \ The category $\modd G_rT$ is a sum of blocks of $\modd (G_r\!\rtimes\!T)$.

{\rm (2)} \ The category $\modd G_rT$ has almost split sequences.

{\rm (3)} \ The canonical restriction functor $\fF : \modd G_rT \lra \modd G_r$ sends indecomposables to indecomposables and almost split sequences to almost split sequences. \hfill 
$\square$ \end{Lemma}

\bigskip
\noindent
Recall from \cite[(II.9.4)]{Ja3} that an object $P \in \modd G_rT$ is projective if and only if it is injective. Thus, $\modd G_rT $ is a Frobenius category (cf.\ \cite[(I.2)]{Ha}), and we can  
speak of  the {\it stable Auslander--Reiten quiver} $\Gamma_s(G_rT)$ of $G_rT$. In view of Lemma \ref{AR-GT2}(3), the functor $\fF$ commutes with the Auslander--Reiten translations 
$\tau_{G_rT}$ and $\tau_{G_r}$ of $\Gamma_s(G_rT)$ and $\Gamma_s(G_r)$, that is,
\[ \tau_{G_r}\circ \fF = \fF \circ \tau_{G_rT}.\]
We record the analogue of \cite[(2.2)]{Fa5}:

\bigskip

\begin{Lemma} \label{AR-GT3} The following statements hold:

{\rm (1)} \ The canonical restriction functor $\fF : \modd G_rT \lra \modd G_r$ induces a morphism $\fF : \Gamma_s(G_rT) \lra \Gamma_s(G_r)$ of stable translation quivers that maps the 
set $x^+$ of successors of an arbitrary vertex $x \in \Gamma_s(G_rT)$ onto the set of successors $\fF(x)^+$ of $\fF(x)$.

{\rm (2)} \ If $\Theta \subseteq \Gamma_s(G_rT)$ is a component, then $\fF(\Theta)$ is a component of $\Gamma_s(G_r)$. \hfill $\square$ \end{Lemma}

\bigskip
\noindent
By work of Gordon and Green \cite[\S1,\S3]{GG2}, the Auslander--Reiten translation $\tau_{G_rT} = \tau_{G_r\rtimes T}|_{\modd G_rT}$ is given by
\[ \tau_{G_rT}(M) = {\rm Tr}(M)^\ast,\]
where ${\rm Tr}(M)$ denotes the {\it transpose} of $M$ (see \cite[(IV.1)]{ARS} for the definition). In view of $\modd G_rT$ being Frobenius, we have natural isomorphisms 
\[{\rm Tr}(M)^\ast \cong \cN(\Omega^2_{G_rT}(M)),\] 
where $\cN = \Hom_{G_r}(-,\Dist(G_r))^\ast$ is the {\it Nakayama functor} of $\modd G_rT$. According to \cite[(I.9.7)]{Ja3}, the group $G$ acts via the character $g \mapsto 
\det(\Ad(g))^{p^r-1}$ on the space of left integrals of $\Dist(G_r)$. The reductivity of $G$ implies that this character is trivial, and we obtain $\cN \cong {\rm id}_{\modd G_rT}$,
cf.\ \cite[(I.8.12)]{Ja3} and \cite{Fa7}. As a result, we have
\[ \tau_{G_rT} (M) \cong \Omega^2_{G_rT}(M)\]
for every $G_rT$-module $M$. 

Replacing $U_0(\fg)$ by $\Dist(G_r)$ in \cite[(2.3)]{Fa5} while observing \cite[(1.3)]{Fa4}, one obtains:

\bigskip  

\begin{Proposition} \label{AR-GT4} Let $\Theta \subseteq \Gamma_s(G_rT)$ be a component. Then the tree class $\bar{T}_\Theta$ of $\Theta$ is a simply laced finite or infinite Dynkin
diagram, a simply laced Euclidean diagram, or $\tilde{A}_{12}$.  \hfill $\square$ \end{Proposition}

\bigskip
\noindent
In contrast to $\modd G_r$, the category $\modd G_rT$ behaves well under passage to coverings. Given a connected smooth reductive group scheme $G$, a connected smooth reductive group 
scheme $\tilde{G}$ is called a {\it covering} of $G$ if there exists a finite diagonalizable subgroup scheme $\tilde{Z} \subseteq \tilde{G}$ such that $\tilde{G}/\tilde{Z} \cong G$ (cf.\ 
\cite[(II.1.17)]{Ja3}).

\bigskip

\begin{Lemma} \label{AR-GT5} Let $\tilde{G}$ be a covering of $G$ with maximal torus $\tilde{T} \subseteq \tilde{G}$ such that $\tilde{T}/\tilde{Z} \cong T$. Then $\modd G_rT$ is the 
sum of those blocks of $\modd \tilde{G}_r\tilde{T}$ whose characters belong to $X(T) \subseteq X(\tilde{T})$. \end{Lemma}

\begin{proof} Owing to \cite[(II.9.7)]{Ja3}, we have $(\tilde{G}_r\tilde{T})/\tilde{Z} \cong G_rT$. By rigidity of tori (cf.\ \cite[(7.7)]{Wa}), the group $\tilde{Z}$ belongs to the center of 
$\tilde{G}$. Thus, $\tilde{Z}$ acts on each module of a block of $\modd \tilde{G}_r\tilde{T}$ via the same character (cf.\ \cite[(II.7.1)]{Ja3}), and our assertion follows. \end{proof} 

\bigskip
\noindent
Recall that an indecomposable module $M$ of a self-injective algebra $\Lambda$ is {\it periodic}, if there exists $m \in \NN$ such that $\Omega^m_\Lambda(M) \cong M$. Periodic
modules have complexity $1$. It thus follows from Lemma \ref{AR-S2} that the simple $G_r$-modules $L_r(\lambda)$ are not periodic. 

The foregoing result enables us to generalize Theorem \ref{AR-S6} to arbitrary reductive groups:

\bigskip 

\begin{Theorem} \label{AR-GT6} Suppose that $G$ is defined over $\FF_p$ and let $\lambda \in X(T)$.

{\rm (1)} \ If the component $\Theta \subseteq \Gamma_s(G_r)$ containing $L_r(\lambda)$ is isomorphic to $\ZZ[A_\infty]$, then $\widehat{L}_r(\lambda)$ is the only simple vertex in 
its AR-component.

{\rm (2)} \ The module $L_r(\lambda)$ is the only simple vertex in its AR-component.\end{Theorem}

\begin{proof} (1) General theory ensures the existence of a covering $\tilde{G}$ such that $(\tilde{G},\tilde{G})$ is simply connected. In view of Lemma \ref{AR-GT5}, we may thus
assume without loss of generality that the derived group $(G,G)$ of $G$ is simply connected. 

Suppose that $\widehat{L}_r(\lambda')$ also belongs to the component $\widehat{\Theta}$ of $\Gamma_s(G_rT)$ containing $\widehat{L}_r(\lambda)$. Thanks to Lemma \ref{AR-GT3}, 
the component $\Theta = \fF(\widehat{\Theta})$ of $\Gamma_s(G_r)$ contains $L_r(\lambda)$ and $L_r(\lambda')$, and a consecutive application of Theorem \ref{AR-S6} and 
\cite[(II.3.10)]{Ja3} implies 
\[ \lambda' = \lambda + p^r\gamma\]
for some $\gamma \in X(T)$. 

In light of Lemma \ref{AR-GT3} and our current assumption, the component $\widehat{\Theta}$ has infinitely many $\tau_{G_rT}$-orbits, and Proposition \ref{AR-GT4} shows that its tree 
class $\bar{T}_{\widehat{\Theta}}$ is isomorphic to $A_\infty, A_\infty^\infty$, or $D_\infty$. Since the number of non-projective direct summands of the middle terms of almost split 
sequences in $\widehat{\Theta}$ takes the values $1$ and $2$, the latter two alternatives cannot occur.

As $L_r(\lambda)$ is quasi-simple, it now follows from Lemma \ref{AR-GT2} that $\widehat{L}_r(\lambda)$ and $\widehat{L}_r(\lambda')$ both lie at ends of $\widehat{\Theta}$.
Consequently, there exists $m \in \ZZ$ with
\[ \Omega_{G_rT}^{2m}(\widehat{L}_r(\lambda)) \cong \widehat{L}_r(\lambda').\]
Application of $\fF$ gives
\[ \Omega_{G_r}^{2m}(L_r(\lambda)) \cong L_r(\lambda') \cong L_r(\lambda).\]
Since the simple $G_r$-modules are not periodic, we conclude $m=0$ and $\lambda = \lambda'$.

(2) Let $\Theta \subseteq \Gamma_s(G_r)$ be the component containing $L_r(\lambda)$. Thanks to Theorem \ref{AR-S3}, we have $\Theta \cong \ZZ[A_\infty], \ZZ[\tilde{A}_{12}]$. In 
the latter case, our assertion follows from Lemma \ref{AR-S5}. 

Suppose that $\Theta \cong \ZZ[A_\infty]$ and denote by $\widehat{\Theta}$ the stable AR-component of $\widehat{L}_r(\lambda)$. The assumption $[L_r(\mu)] \in \Theta = 
\fF(\widehat{\Theta})$ implies the existence of $\gamma \in X(T)$ with $[\widehat{L}_r(\mu+p^r\gamma)] \in \widehat{\Theta}$ (see \cite[(II.3.10)]{Ja3}). According to (1), we have 
$\lambda = \mu +p^r\gamma$, whence $L_r(\lambda) \cong L_r(\mu)$, as desired.  \end{proof} 

\bigskip

\begin{Remark} If $\Theta \cong \ZZ[\tilde{A}_{12}]$, then the AR-component containing $\widehat{L}_r(\lambda)$ may have several simple vertices, see \cite[p.749f]{Fa5}.
\end{Remark}

\bigskip
\noindent
We turn to the consideration of $\Delta$-good modules of the highest weight category $\modd G_rT$. Using the rank varieties $V_r(G)_M$, the arguments employed in \cite[(4.1)]{Fa5} 
imply:

\bigskip

\begin{Theorem} \label{AR-GT7} Let $M$ be an indecomposable $G_rT$-module, $\Theta \subseteq \Gamma_s(G_rT)$ and $\Psi \subseteq \Gamma_s(G_r)$ be the stable AR-components 
containing $[M]$ and $[\fF(M)]$, respectively.

{\rm (1)} \ Every vertex of $\Psi$ has a $G_rT$-structure.

{\rm (2)} \ If $M$ affords a $\widehat{Z}_r$-filtration, so does every indecomposable $G_rT$-module belonging to $\Theta$.

{\rm (3)} \ If $\fF(M)$ affords a $Z_r$-filtration, so does every indecomposable $G_r$-module belonging to $\Psi$. \hfill  $\square$ \end{Theorem}

\bigskip
\noindent
We conclude this section by recording two elementary, but useful observations:

\bigskip

\begin{Lemma} \label{AR-GT8} Let
\[ (0) \lra M' \lra M \lra M'' \lra (0) \]
be an exact sequence of $G_rT$-modules. If two modules of the sequence belong to $\cF(\Delta)$, so does the third. In that case, we have
\[ [M\!:\!\widehat{Z}_r(\lambda)] = [M'\!:\!\widehat{Z}_r(\lambda)] + [M''\!:\!\widehat{Z}_r(\lambda)]\]
for every $\lambda \in X(T)$.\end{Lemma}

\begin{proof} Recall that $\modd B_r^-$ is a Frobenius category. If two modules of the above sequence are $B^-_r$-injective, so is the third. Thus, the first assertion follows directly from
\cite[(II.11.2)]{Ja3}. By the same token, we have
\[ [M\!:\!\widehat{Z}_r(\lambda)] = \dim_k\Hom_{G_rT}(M,\widehat{Z}'_r(\lambda)),\]
so that the additivity of the filtration multiplicities follows by applying $\Hom_{G_rT}(-,\widehat{Z}'_r(\lambda))$, while observing 
$\Ext^1_{G_rT}(\widehat{Z}_r(\mu)),\widehat{Z}'_r(\lambda)) = (0)$ for all $\mu \in X(T)$, see \cite[(II.9.9)]{Ja3}. \end{proof}

\bigskip

\begin{Lemma} \label{AR-GT9} Let $M \in \cF(\Delta)$. Then $\Top_{G_rT}(M)$ is a direct summand of 
\[ \bigoplus_{\lambda \in \supp_\Delta(M)} [M\!:\!\widehat{Z}_r(\lambda)] \widehat{L}_r(\lambda).\] \end{Lemma}

\begin{proof} We use induction on the  $\widehat{Z}_r$-filtration length $\ell_\Delta(M)$ of $M$, the case $\ell_\Delta(M) = 1$ being a consequence of \cite[(II.9.6(4))]{Ja3}. An exact 
sequence
\[ (0) \lra M' \lra M \lra \widehat{Z}_r(\mu) \lra (0)\]
of $G_rT$-modules induces an exact sequence
\[ \Top_{G_rT}(M') \lra \Top_{G_rT}(M) \lra \widehat{L}_r(\mu) \lra (0),\]
so that $\Top_{G_rT}(M)$ is a direct summand of $\Top_{G_rT}(M')\oplus \widehat{L}_r(\mu)$. Lemma \ref{AR-GT8} yields
\[ [M\!:\!\widehat{Z}_r(\lambda)] = [M'\!:\!\widehat{Z}_r(\lambda)]+ \delta_{\lambda,\mu},\]
and by inductive hypothesis, $\Top_{G_rT}(M)$ is a direct summand of
\[  \bigoplus_{\lambda \in \supp_\Delta(M')} [M'\!:\!\widehat{Z}_r(\lambda)] \widehat{L}_r(\lambda) \oplus \widehat{L}_r(\mu)
=  \bigoplus_{\lambda \in \supp_\Delta(M)} [M\!:\!\widehat{Z}_r(\lambda)] \widehat{L}_r(\lambda),\]
as desired. \end{proof}

\bigskip

\section{Auslander-Reiten Components of Verma Modules} \label{s:AR-Ve}
Given a component $\Theta \subseteq \Gamma_s(G_r)$, we recall that any two vertices of $\Theta$ have the same rank variety $V_r(G)_\Theta$ (cf.\ \cite[(1.1)]{Fa4}). The quiver 
$\Gamma_s(G_r)$ is known to satisfy an analogue of Webb's Theorem \cite{We} for finite groups, that is, the tree classes of the components of $\Gamma_s(G_r)$ are finite or infinite Dynkin 
diagrams, or Euclidean diagrams (cf.\ \cite[(1.3)]{Fa4}). Our next lemma further reduces the possibilities for the components $\Theta_r(\lambda)$ containing the isoclasses $[Z_r(\lambda)]$.

The group $G$ acts on $G_r$ via conjugation. This action induces operations on $V_r(G)$ and $\Dist(G_r)$, where we denote the automorphism of $\Dist(G_r)$ corresponding to $g$ by 
$\Ad(g)$. Letting $M^{(g)}$ be the pull-back of a $G_r$-module $M$ along $\Ad(g)^{-1}$, we obtain actions of $G$ on $\modd G_r$ and $\Gamma_s(G_r)$ via auto-equivalences and 
isomorphisms of stable representation quivers, respectively. Note that
\[ V_r(G)_{M^{(g)}} = g\dact V_r(G)_{M} \ \ \ \ \forall \ g \in G.\]
Since $G$ is connected, the isoclasses of the simple $G_r$-modules are fixed by the above action, so that the rank varieties of these modules are $G$-stable subsets of $V_r(G)$. 

A component of $\Gamma_s(G_r)$ is {\it regular} if it is not attached to a principal indecomposable module.   

\bigskip

\begin{Lemma} \label{AR-V1} Let $\lambda \in X(T)$ be a weight such that $Z_r(\lambda)$ is not projective. Then the following statements hold:

{\rm (1)} \ The component $\Theta_r(\lambda)$ is regular.

{\rm (2)} \ The isoclass $[Z_r(\lambda)]$ is a vertex of $\Theta_r(\lambda)$ of minimal dimension.

{\rm (3)} \ If $\dim V_r(G)_{\Theta_r(\lambda)} = 1$, then $\Theta_r(\lambda) \cong \ZZ[A_\infty]/(\tau^{p^n})$ for some $n \in \{0,\ldots,r-1\}$.

{\rm (4)} \ If $\dim V_r(G)_{\Theta_r(\lambda)} \ge 2$, then $\Theta_r(\lambda) \cong \ZZ[A_\infty], \ZZ[D_\infty]$, or $\ZZ[A_\infty^\infty]$. \end{Lemma}

\begin{proof} Thanks to \cite[(1.5,1.14)]{SFB1} and \cite[(5.2)]{SFB2}, there is a commutative diagram 
\[ \begin{CD} V_r(U) @> \Psi_U>> \cV_{U_r}(k) \\
        @V\iota VV @V {\rm res}VV\\
 V_r(G) @>\Psi_G>> \cV_{G_r}(k), \end{CD}\]
 where $\iota$ is injective and the horizontal arrows are homeomorphisms. Owing to Proposition \ref{Pr2} we have  $\cV_{G_r}(Z_r(\lambda)) \subseteq \im {\rm res}$. Thus, 
 \cite[(6.8)]{SFB2} gives $V_r(G)_{\Theta_r(\lambda)} \subseteq V_r(U)$. 

Let $g_0 \in G$ be a representative of the longest element $w_0$ of the Weyl group $W$ of $G$, so that $g_0Ug_0^{-1} = U^-$, cf.\ \cite[(II.1.8)]{Ja3}. For any non-trivial element 
$\varphi \in V_r(G)_{\Theta_r(\lambda)} \subseteq V_r(U)$, we have $g_0\dact \varphi \in V_r(U^-)$, so that $g_0\dact \varphi \not \in V_r(G)_{\Theta_r(\lambda)}$. Thus,
$V_r(G)_{\Theta_r(\lambda)}$ is not $G$-invariant.

(1) Since rank varieties of non-regular components are necessarily $G$-invariant (cf.\ \cite[\S7]{Fa3}), assertion (1) follows.

(2) This is a direct consequence of Theorem \ref{AR-GT7}(3).

(3) A consecutive application of \cite[(2.1)]{Fa4} and \cite[(4.1)]{Fa4} shows that $\Theta_r(\lambda) \cong  \ZZ[A_\infty]/(\tau^m)$ for some $m \in \NN$. According to Lemma \ref{AR-S1}, $\Dist(G_r)$ is a symmetric algebra, so that $\tau_{G_r} = \Omega^2_{G_r}$.  The Friedlander-Suslin Theorem \cite[(1.5)]{FS} shows that the even cohomology ring
$\HH^\bullet(G_r,k)$ is a finitely generated module of a graded subalgebra, whose generators have degree $2p^i$ for $0\le i \le r\!-\!1$. Thanks to \cite[(5.10.6)]{Be2}, the periods of
the periodic modules are divisors of these numbers. This implies our assertion concerning $n$.
 
(4) Suppose that assertion (4) does not hold. The classification of AR-components (cf.\ \cite[(4.1)]{Fa4}) now implies that $\Theta_r(\lambda)$ is a Euclidean component. As such, it is not regular, contradicting part (1). \end{proof}

\bigskip
\noindent
Our next result implies that the vertex $[\widehat{Z}_r(\lambda)]$ has exactly one predecessor in $\Gamma_s(G_rT)$, thereby showing in particular that $\Theta_r(\lambda) \not \cong 
\ZZ[A_\infty^\infty]$. 

\bigskip

\begin{Theorem} \label{AR-V2} Suppose that $\widehat{Z}_r(\lambda)$ is not projective and let $\widehat{E}$ be the middle term of the almost split sequence
\[ (0) \lra \Omega^2_{G_rT}(\widehat{Z}_r(\lambda)) \lra \widehat{E} \lra \widehat{Z}_r(\lambda) \lra (0).\]
Then the following statements hold:

{\rm (1)} \ $\supp_\Delta(\widehat{E}) \subseteq X(T)_{\ge \lambda}$ and $[\widehat{E}\!:\!\widehat{Z}_r(\lambda)] = 1$.

{\rm (2)} \ The module $\widehat{E}$ is indecomposable. \end{Theorem}

\begin{proof} We consider a minimal projective presentation 
\[ (0) \lra \Omega^2_{G_rT}(\widehat{Z}_r(\lambda)) \lra \widehat{P}_1 \lra \widehat{P}_r(\lambda) \lra \widehat{Z}_r(\lambda) \lra (0)\]
of the $G_rT$-module $\widehat{Z}_r(\lambda)$. Since $[\widehat{Z}_r(\lambda)\!:\!\widehat{L}_r(\lambda)] = 1$ (cf.\ \cite[(II.9.6)]{Ja3}), Lemma \ref{AR-GT1} in conjunction with 
Lemma \ref{AR-GT8} implies the inclusion $\supp_\Delta(\Omega_{G_rT}(\widehat{Z}_r(\lambda))) \subseteq X(T)_{>\lambda}$. Owing to Lemma \ref{AR-GT9}, the module 
$\Top_{G_rT}(\Omega_{G_rT}(\widehat{Z}_r(\lambda)))$ is a direct summand of 
\[  \bigoplus_{\mu \in \supp_\Delta(\Omega_{G_rT}(\widehat{Z}_r(\lambda)))} [\Omega_{G_rT}(\widehat{Z}_r(\lambda))\!:\!\widehat{Z}_r(\mu)] \widehat{L}_r(\mu).\] 
Another application of Lemma \ref{AR-GT1} therefore yields $\supp_\Delta(\widehat{P}_1) \subseteq X(T)_{>\lambda}$. In particular, the module $\widehat{Z}_r(\lambda)$ is not a 
filtration factor of $ \Omega^2_{G_rT}(\widehat{Z}_r(\lambda))$ and Lemma \ref{AR-GT8} implies
\[ [\widehat{E}\!:\!\widehat{Z}_r(\lambda)] = 1.\]
In virtue of Theorem \ref{AR-GT7}(2), every indecomposable constituent of $\widehat{E}$ is $\widehat{Z}_r$-filtered, so that each irreducible morphism terminating in 
$\widehat{Z}_r(\lambda)$ is surjective (see \cite[(V.5.1)]{ARS}). Thus, if $\widehat{E}_i$ is an indecomposable constituent of $\widehat{E}$, then, by \cite[(V.5.3)]{ARS}, our 
almost split sequence induces an exact sequence
\[ (0) \lra X \lra \widehat{E}_i \lra \widehat{Z}_r(\lambda) \lra (0),\]
and Lemma \ref{AR-GT8} gives rise to
\[ [\widehat{E}_i\!:\!\widehat{Z}_r(\lambda)] = [X\!:\!\widehat{Z}_r(\lambda)] + [\widehat{Z}_r(\lambda)\!:\!\widehat{Z}_r(\lambda)] = [X\!:\!\widehat{Z}_r(\lambda)] + 1.\]
Consequently, $\widehat{E}$ can only have one indecomposable constituent, so that (2) follows. 

A twofold application of Lemma \ref{AR-GT8} yields
\[ \supp_\Delta(\widehat{E}) = \supp_\Delta(\widehat{Z}_r(\lambda))\cup \supp_\Delta(\Omega_{G_rT}^2(\widehat{Z}_r(\lambda))) \subseteq
\supp_\Delta(\widehat{Z}_r(\lambda))\cup \supp_\Delta(\widehat{P}_1) \subseteq X(T)_{\ge \lambda},\]
as desired. \end{proof}

\bigskip

\begin{Corollary} \label{AR-V3} Suppose that $Z_r(\lambda)$ is not projective. The middle term $E$ of the almost split sequence
\[ (0) \lra \Omega^2_{G_r}(Z_r(\lambda)) \lra E \lra Z_r(\lambda) \lra (0)\]
is indecomposable. \end{Corollary}

\begin{proof} According to \cite[(II.9.4)]{Ja3} the $G_rT$-module $\widehat{Z}_r(\lambda)$ is not projective. We consider the almost split sequence
\[(0) \lra \Omega^2_{G_rT}(\widehat{Z}_r(\lambda)) \lra \widehat{E} \lra \widehat{Z}_r(\lambda) \lra (0)\]
terminating in $\widehat{Z}_r(\lambda)$. Since $Z_r(\lambda) = \fF(\widehat{Z}_r(\lambda))$, Lemma \ref{AR-GT2}(3) shows that
\[(0) \lra \fF(\Omega^2_{G_rT}(\widehat{Z}_r(\lambda))) \lra \fF(\widehat{E}) \lra Z_r(\lambda) \lra (0)\]
is the almost split sequence terminating in $Z_r(\lambda)$. In view of Theorem \ref{AR-V2}, Lemma \ref{AR-GT2}(3) and \cite[(V.1.16)]{ARS}, its middle term $E \cong \fF(\widehat{E})$ 
is indecomposable. \end{proof}

\bigskip
\noindent

\begin{Theorem} \label{AR-V4} Let $\lambda \in X(T)$. The baby Verma module $Z_r(\lambda)$ is either projective or quasi-simple. \end{Theorem}

\begin{proof} In view of Lemma \ref{AR-V1} and Corollary \ref{AR-V3}, it remains to rule out the case where $\Theta_r(\lambda) \cong \ZZ[D_\infty]$.

Let $\widehat{\Theta}_r(\lambda) \subseteq \Gamma_s(G_rT)$ be the component containing the vertex $[\widehat{Z}_r(\lambda)]$. Thanks to Lemma \ref{AR-GT3}(2),  we have 
$\Theta_r(\lambda) = \fF(\widehat{\Theta}_r(\lambda))$.  Assuming $\Theta_r(\lambda) \cong \ZZ[D_\infty]$, we consider the almost split sequence
\[ (\ast) \ \ \ \ \ \ \ \ \ (0) \lra \Omega^2_{G_rT}(\widehat{Z}_r(\lambda)) \lra \widehat{E} \lra \widehat{Z}_r(\lambda) \lra (0)\]
terminating in $\widehat{Z}_r(\lambda)$. By assumption, there exists an indecomposable $G_r$-module $Y_r(\lambda) \not \cong Z_r(\lambda)$ such that $[E] := [\fF(\widehat{E})]$ is 
the only predecessor of $[Y_r(\lambda)] \in \Theta_r(\lambda)$. Let $[\widehat{Y}_r(\lambda)] \in \widehat{\Theta}_r(\lambda)$ be a pre-image of $[Y_r(\lambda)]$ under $\fF$, and 
consider the almost split sequence
\[ (\ast\ast) \ \ \ \ \ \ \ \ \ (0) \lra \Omega^2_{G_rT}(\widehat{Y}_r(\lambda)) \lra \widehat{D} \lra \widehat{Y}_r(\lambda) \lra (0)\] 
terminating in $\widehat{Y}_r(\lambda)$. 

We now proceed in several steps:

(i) \ {\it We have $\widehat{D} \cong \widehat{E}$}.

\smallskip
\noindent
Lemma \ref{AR-GT2}(3) provides isomorphisms $\fF(\widehat{D}) \cong E \cong \fF(\widehat{E})$, and \cite[(4.1)]{GG1} implies the existence of $\mu \in X(T)$ such that $\widehat{D} 
\cong \widehat{E}\!\otimes_k\!k_\mu$. 

Since $\fF$ sends the $\tau_{G_rT}$-orbits in $\widehat{\Theta}_r(\lambda)$ onto the $\tau_{G_r}$-orbits in $\Theta_r(\lambda)$ (see Lemma \ref{AR-GT2}), we conclude that the 
component $\widehat{\Theta}_r(\lambda)$ has infinitely many $\tau_{G_rT}$-orbits. Owing to Proposition \ref{AR-GT4}, the tree class $\bar{T}_{\widehat{\Theta}_r(\lambda)}$ is 
therefore isomorphic to $A_\infty, \, A_\infty^\infty$, or $D_\infty$. However, in components with the former two tree types, every vertex has at most two successors. Thanks to Lemma 
\ref{AR-GT3}, the component $\Theta_r(\lambda)$ would also enjoy this property. As $\Theta_r(\lambda) \cong \ZZ[D_\infty]$, this is not possible, so that $\widehat{\Theta}_r(\lambda)$ 
also has tree class $D_\infty$. 

By general theory, the auto-equivalence $[M] \mapsto [M\!\otimes_k\!k_\mu]$ induces an automorphism of $D_\infty$. As this map fixes the node corresponding to $\widehat{E}$,
there exists $m_{\widehat{E}} \in \ZZ$ such that $\tau^{m_{\widehat{E}}}_{G_rT}(\widehat{E}) \cong \widehat{E}\!\otimes _k\!k_\mu$. Application of $\fF$ gives
\[ \Omega_{G_r}^{2\,m_{\widehat{E}}}(E) \cong \fF( \widehat{E}\!\otimes _k\! k_\mu) \cong \fF(\widehat{D}) \cong E.\]
The assumption $m_{\widehat{E}} \ne 0$ thus yields the periodicity of $E$, whence
\[ 1 = \dim V_r(G)_E = \dim V_r(G)_{\Theta_r(\lambda)} = \dim V_r(G)_{Z_r(\lambda)}.\]
We may now apply Lemma \ref{AR-V1} to reach a contradiction. Consequently, $m_{\widehat{E}}=0$, whence $\widehat{E} \cong \widehat{E}\!\otimes_k\!k_\mu \cong \widehat{D}$. 
\hfill $\Diamond$ 

\medskip

(ii) \ {\it We have $[\widehat{Y}_r(\lambda)\!:\!\widehat{Z}_r(\lambda)] = 1$.}

\smallskip
\noindent
Consider the almost split sequence
\[  (0) \lra \Omega^2_{G_rT}(\widehat{Y}_r(\lambda)) \lra \widehat{E} \lra \widehat{Y}_r(\lambda) \lra (0).\] 
By virtue of Lemma \ref{AR-GT8} and Theorem \ref{AR-V2}(1), the assumption $(\widehat{Y}_r(\lambda)\!:\!\widehat{Z}_r(\lambda)) = 0$ implies 
\[\supp_\Delta(\widehat{Y}_r(\lambda))\subseteq X(T)_{>\lambda}.\] 
Thus, by Lemma \ref{AR-GT9}, the highest weights $\mu$ of the constituents $\widehat{L}_r(\mu)$ of $\Top(\widehat{Y}_r(\lambda))$ belong to $X(T)_{>\lambda}$. Owing to Lemma 
\ref{AR-GT1}, the projective cover $\widehat{P}(\widehat{Y}_r(\lambda))$ of $\widehat{Y}_r(\lambda)$ satisfies $\supp_\Delta (\widehat{P}(\widehat{Y}_r(\lambda)))
\subseteq X(T)_{>\lambda}$ and, thanks to Lemma \ref{AR-GT8}, the submodule $\Omega_{G_rT}(\widehat{Y}_r(\lambda))$ inherits this property. A repetition of this argument implies
$\supp_\Delta (\Omega^2_{G_rT}(\widehat{Y}_r(\lambda))) \subseteq X(T)_{>\lambda}$. In virtue of Lemma \ref{AR-GT8}, our almost split sequence then yields
\[[\widehat{E}\!:\!\widehat{Z}_r(\lambda)] = [\Omega^2_{G_rT}(\widehat{Y}_r(\lambda))\!:\!\widehat{Z}_r(\lambda)] + [\widehat{Y}_r(\lambda)\!:\!\widehat{Z}_r(\lambda)] = 
0,\]
which contradicts Theorem \ref{AR-V2}(1). Consequently, $[\widehat{Y}_r(\lambda)\!:\!\widehat{Z}_r(\lambda)] \ge 1$, and the reverse inequality follows from Lemma \ref{AR-GT8} and 
Theorem \ref{AR-V2}(1). \hfill $\Diamond$ 

\medskip
(iii) \ {\it There exists a surjection $g : \widehat{Y}_r(\lambda) \lra \widehat{Z}_r(\lambda)$}. 

\smallskip
\noindent
In virtue of Theorem \ref{AR-V2} and (ii), we have $\supp_\Delta (\widehat{Y}_r(\lambda)) \subseteq X(T)_{\ge \lambda}$ as well as $[\widehat{Y}_r(\lambda)\!:\!
\widehat{Z}_r(\lambda)] = 1$. Consequently, there exist submodules $N \subseteq M \subseteq \widehat{Y}_r(\lambda)$ such that

(a) \ $M/N \cong \widehat{Z}_r(\lambda)$, and

(b) \ $\widehat{Y}_r(\lambda)/M$ belongs to $\cF(\Delta)$ with $\supp_\Delta(\widehat{Y}_r(\lambda)/M) \subseteq X(T)_{>\lambda}$.

\noindent
There results an exact sequence
\[ (0) \lra \widehat{Z}_r(\lambda) \lra \widehat{Y}_r(\lambda)/N \lra \widehat{Y}_r(\lambda)/M \lra (0).\]
Lemma \ref{AR-GT1} implies $\Ext^1_{G_rT}(\widehat{Z}_r(\mu),\widehat{Z}_r(\lambda)) = (0)$ for all $\mu \in X(T)_{>\lambda}$. Hence our sequence splits, and
\[ \widehat{Y}_r(\lambda)/N \cong (\widehat{Y}_r(\lambda)/M) \oplus \widehat{Z}_r(\lambda).\]
We may now define $g := {\rm pr} \circ \pi$ to be the composite of the projection ${\rm pr} : \widehat{Y}_r(\lambda)/N \lra \widehat{Z}_r(\lambda)$ with the canonical map $\pi : 
\widehat{Y}_r(\lambda) \lra \widehat{Y}_r(\lambda)/N$. \hfill $\Diamond$

\medskip

(iv) \ {\it We have $\dim_k \Hom_{G_rT}(\widehat{E},\widehat{Z}_r(\lambda)) = 1$}.

\smallskip
\noindent
We write $X := \Omega^2_{G_rT}(\widehat{Z}_r(\lambda))$ for notational convenience. Application of the left exact functor $\Hom_{G_rT}(-,\widehat{Z}_r(\lambda))$ to the almost split 
sequence ($\ast$) terminating in $\widehat{Z}_r(\lambda)$ gives an exact sequence
\[ (0) \lra \Hom_{G_rT}(\widehat{Z}_r(\lambda),\widehat{Z}_r(\lambda)) \lra  \Hom_{G_rT}(\widehat{E},\widehat{Z}_r(\lambda)) \lra \Hom_{G_rT}(X,\widehat{Z}_r(\lambda)).\]
Since $\Top_{G_rT}(X)$ is a direct summand of $\bigoplus_{\mu > \lambda} [X\!:\!\widehat{Z}_r(\mu)] \widehat{L}_r(\mu)$, and all weights of $\widehat{Z}_r(\lambda)$ belong to 
$X(T)_{\le \lambda}$, it follows that any homomorphism $X \lra \widehat{Z}_r(\lambda)$ sends the generators of $X$ to zero. Consequently, $\Hom_{G_rT}(X,\widehat{Z}_r(\lambda)) = 
(0)$, implying
\[ \Hom_{G_rT}(\widehat{E},\widehat{Z}_r(\lambda)) \cong \Hom_{G_rT}(\widehat{Z}_r(\lambda),\widehat{Z}_r(\lambda)) \cong k,\]
as desired. \hfill $\Diamond$

\medskip
\noindent
Let $f : \widehat{E} \lra \widehat{Z}_r(\lambda)$ and $h : \widehat{E} \lra \widehat{Y}_r(\lambda)$ be the irreducible morphisms given by the almost split sequences ($\ast$)  and ($\ast\ast$), terminating in $\widehat{Z}_r(\lambda)$ and $\widehat{Y}_r(\lambda)$, respectively. By virtue of (iii) and (iv), there exists $\alpha \in k \setminus \{0\}$ such that
\[ g \circ h = \alpha \, f.\] 
Since $\alpha \, f$ is irreducible and $h$ is not split injective, the map $g$ is split surjective, \cite[(V.5)]{ARS}. Consequently, $\widehat{Z}_r(\lambda) \cong
\widehat{Y}_r(\lambda)$, a contradiction. \end{proof}

\bigskip

\begin{Corollary} \label{AR-V5} Given $\lambda \in X(T)$, the module $\widehat{Z}_r(\lambda)$ is either projective or quasi-simple. \end{Corollary}

\begin{proof} Assume that $\widehat{Z}_r(\lambda)$ is not projective. Thanks to Theorem \ref{AR-V2}, it suffices to show that the tree class $\widehat{T}_r(\lambda)$ of the component $\widehat{\Theta}_r(\lambda) \subseteq \Gamma_s(G_rT)$ containing $[\widehat{Z}_r(\lambda)]$ is isomorphic to $A_\infty$.

By Theorem \ref{AR-V4}, the component $\Theta_r(\lambda)$ has tree class $A_\infty$ and thus has infinitely many $\tau_{G_r}$-orbits. In view of Lemma \ref{AR-GT3}, the component $\widehat{\Theta}_r(\lambda)$ also has infinitely many $\tau_{G_rT}$-orbits, and an application of Proposition \ref{AR-GT4} now implies $\widehat{T}_r(\lambda) = A_\infty, A_\infty^\infty$, or $D_\infty$. By Theorem \ref{AR-V2} and the proof of Theorem \ref{AR-V4}, the latter two alternatives cannot occur. 
\end{proof}

\bigskip

\begin{Remark} Aside from simple modules and baby Verma modules, the so-called {\it Weyl modules} $V(\lambda)$ (see \cite[(II.2.13)]{Ja3}) and their duals form another important class 
of $G_r$-modules. For the group scheme $\SL(2)_1$, these modules belong to components of type $\ZZ[\tilde{A}_{12}]$, cf. \cite[(4.3.1)]{Fa6}. Since the support varieties of Weyl modules
are $G$-invariant, they will in general have dimension $\ge 3$, so that \cite[(2.2)]{Fa4} ensures that Weyl modules usually belong to components of type $\ZZ[A_\infty]$. However, their position within these components is not known. \end{Remark} 

\bigskip

\section{Varieties for Verma Modules} \label{s:grvar}
We begin by generalizing \cite[(2.3)]{FR} to Verma modules of higher Frobenius kernels. Recall that $h$ is the Coxeter number of the reductive group $G$. Unless mentioned otherwise, the 
prime $p$ is assumed to be good for $G$.

The interested reader may compare the following result with \cite[(4.1.2)]{NPV}. Under the additional hypothesis of $\cV_{U_r}(k)$ being irreducible, the implication (3) $\Rightarrow$ (1) provides the converse to part (b) of that result. Also note that the assumption $p\ge h$ implies that $\cV_{U_1}(k)$ is irreducible.

\bigskip

\begin{Lemma} \label{VV1}  Suppose $p \ge h$ and that $\cV_{U_r}(k)$ is irreducible. Let $r \le s$ and  $\lambda \in X(T)$ be a weight. Then the following statements are equivalent:

{\rm (1)} \ $\cV_{G_r}(Z_s(\lambda)) = \cV_{U_r}(k)$.

{\rm(2)} \ $\cx_{G_r}(Z_s(\lambda)) = \dim \cV_{U_r}(k)$.

{\rm (3)} \ $\lambda$ is $p$-regular. \end{Lemma}

\begin{proof} (1) $\Rightarrow$ (2). Thanks to \cite[(1.5.2)]{FS} the cohomology ring $\HH^{\rm ev}(G_r, k)$ is finitely generated.  Thus, the arguments of \cite[(2.3)]{Ca1} show that 
the complexity of the $G_r$-module $Z_s(\lambda)$ coincides with the dimension of $\cV_{G_r}(Z_s(\lambda))$.

(2) $\Rightarrow$ (3). The $G_rT$-module $\widehat{Z}_s(\lambda)$ is $B_r^-$-projective, and \cite[(II.11.2)]{Ja3} provides a $\widehat{Z}_r$-filtration of 
$\widehat{Z}_s(\lambda)|_{G_rT}$. Applying the exact functor $\fF$, we conclude that $Z_s(\lambda)|_{G_r}$ has a $Z_r$-filtration, so that Proposition \ref{Pr2} gives 
\[ \cV_{G_r}(Z_s(\lambda)) \subseteq \cV_{U_r}(k) \]
for any $\lambda \in X(T)$. Since the variety $\cV_{U_r}(k)$ is irreducible, (2) implies that the above inclusion is in fact an equality. In view of \cite[(4.1.2)]{NPV}, we obtain the 
$p$-regularity of $\lambda$.

(3) $\Rightarrow$ (1). Following the proof of \cite[(2.3)]{FR}, we first consider the case where $\lambda = 0$. The augmentation map $\varepsilon : \Dist(G_r) \lra k$ furnishes a split surjective homomorphism
\[ f : Z_s(0) \lra k \ \ ; \ \  u\otimes \alpha \mapsto \varepsilon(u)\alpha\]
of $\Dist(U_r)$-modules. Thus, $k$ is a direct summand of $Z_s(0)|_{U_r}$, so that $\cV_{U_r}(k) \subseteq \cV_{U_r}(Z_s(0))$. The foregoing observations now yield
\[ \cV_{G_r}(Z_s(0)) = {\cV}_{U_r}(k).\]
To transfer this result to other $p$-regular weights, we shall use translation functors $T_\lambda^\mu$, cf.\ \cite[(II.9.22)]{Ja3}. Their definition in conjuction with \cite[(7.2)]{SFB2}, 
yields 
\[\cV_{G_r}(T_\lambda^\mu(M)) \subseteq \cV_{G_r}(M).\]
Now let $\lambda$ be another weight in the standard alcove 
\[ C_0 := \{ \lambda \in X(T)\!\otimes_\ZZ\!\RR \mid 0 < \langle \lambda + \rho,\alpha^\vee \rangle < p \ \ \ \ \forall \ \alpha \in \Psi^+\}\]
(such weights are necessarily $p$-regular). The translation functor $T_0^\lambda$ provides a categorical equivalence between the blocks of $\Dist(G_r)$ associated to $\lambda$ and $0$ 
and sends $Z_s(0)$ onto $Z_s(\lambda)$ (cf.\ \cite[(II.9.22(1)-(3))]{Ja3}). Consequently, the $G_r$-support varieties of $Z_s(\lambda)$ and $Z_s(0)$ are equal. 

By the same token, $\cV_{G_r}(Z_s(\lambda))$ and $\cV_{G_r}(Z_s(\gamma))$ coincide, whenever $\lambda$ and $\gamma$ are $p$-regular elements belonging to the same alcove.

Let $C$ be an alcove such that $\cV_{G_r}(Z_s(\lambda)) =  \cV_{U_r}(k)$ for some (and hence every) $p$-regular element $\lambda \in C$. We choose an alcove $D$ adjacent to $C$ and 
let $s_F$ be the unique affine reflection relative to the adjoining wall $\overline F =\overline C\cap \overline D$ (facet of codimension $1$). Note that $s_F\dact \lambda \in D$ is a 
$p$-regular element. Next, we choose an element $w \in W_p$ such that $w\dact C = C_0$. Then $w\dact D$ is adjacent to $C_0$ with adjoining wall $w\dact F$ and associated reflection 
$s' = ws_Fw^{-1}$. Observe that $w^{-1}\dact s'\dact (w\dact \lambda)=  s_F\dact \lambda$. Now let $\mu \in w\dact F$ be a weight on the wall $w\dact F$.  Suppose that $\lambda <  
s_F\dact \lambda$. An application of \cite[(II.9.22(3))]{Ja3} provides an exact sequence
\[ (0) \lra Z_s(\lambda) \lra T_\mu^{w\dact\lambda}(Z_s(w^{-1}\dact\mu)) \lra Z_s(s_F\dact \lambda) \lra (0)\]
of $G_s$-modules, so that
\[\cV_{G_r}(Z_s(\lambda)) \subseteq \cV_{G_r}(T_\mu^{w\dact\lambda}(Z_s(w^{-1}\dact\mu))) \cup \cV_{G_r}(Z_s(s_F \dact \lambda))\subseteq \cV_{G_r}(Z_s(w^{-1}\dact\mu)) 
\cup \cV_{G_r}(Z_s(s_F \dact \lambda)). \]
Since $\cV_{G_r}(Z_s(\lambda)) = \cV_{U_r}(k)$ is irreducible and $w^{-1}\dact\mu \in w\dact F$ is not $p$-regular, we have $\cV_{G_r}(Z_s(\lambda)) = \cV_{G_r}(Z_s(s_F \dact
\lambda))$, cf.\ \cite[(4.1.2)]{NPV}. Therefore, $\cV_{G_r}(Z_s(\gamma))=  \cV_{U_r}(k)$ for any $p$-regular element $\gamma \in D$. The case $\lambda > s_F\dact \lambda$ can be 
treated similarly.

Since $W_p$ acts transitively on the set of alcoves, we are done. \end{proof}

\bigskip
\noindent
We shall use the foregoing result mainly for Levi subgroups defined by simple roots. The following example provides another case, where $\cV_{U_r}(k)$ is irreducible.

\bigskip

\begin{Example} Suppose that $p \ge 3$ and consider the group $G = \SL(3)$. We let $U$ denote the unipotent subgroup of upper triangular matrices. Then $\fu := \Lie(U)$
is isomorphic to the $3$-dimenisonal $p$-unipotent Heisenberg algebra with trivial $p$-map and multiplication relative to the standard basis $\{a,b,c\}$ defined via
\[ [a,b] = c \ \ \text{and} \ \ [a,c] = 0 = [b,c].\]
For $r\ge 2$ and any subspace $\fv \subseteq \fu$, we consider the commuting variety
\[ \cC^r(\fv) := \{(v_0,\ldots,v_{r-1}) \in \fv^r \ \mid \ [v_i,v_j] = 0 \ \text{for} \ i \ne j\}.\]
Writing $v_j := x_ja+y_jb+z_jc$, we obtain
\[(\ast) \ \ \ \ \ \ \ \ \ \ \ [v_i,v_j] = (x_iy_j-x_jy_i)c\]
Thus, for $\fv := ka\oplus kb$, the variety $\cC^r(\fv)$ is defined by the determinantal ideal 
\[ I_2(A) \subseteq k[X_0,\ldots,X_{r-1},Y_0,\ldots,Y_{r-1}],\] 
given by the $(2\!\times\! 2)$-minors of the $(r\!\times\!r)$-matrix
\[ A := \left(\begin{array}{ccc} X_0 & \ldots & X_{r-1} \\ Y_0 & \ldots & Y_{r-1}\end{array}\right).\]
Thanks to \cite[(7.3.1)]{BH}, the determinantal variety $C^r(\fv)$ is irreducible.

In view of ($\ast$), the isomorphism
\[ \fv^r \times (kc)^r \lra \fu^r \ \ ; \ \ (v_0,\ldots,v_{r-1},w_0,\ldots,w_{r-1}) \mapsto (v_0+w_0,\ldots, v_{r-1}+w_{r-1})\]
maps the irreducible variety $\cC^r(\fv)\times (kc)^r$ onto $\cC^r(\fu)$, so that $\cC^r(\fu)$ is also irreducible. By applying \cite[(1.7),(1.8)]{SFB1} and \cite[(5.2)]{SFB2}
successively, we conclude that the support variety $\cV_{U_r}(k)$ is irreducible. Owing to \cite[(7.3.1)]{BH}, we also have $\dim \cV_{U_r}(k) = 2r+1$.\end{Example}

\bigskip
\noindent
Given a root $\alpha \in \Psi$, we let $U_\alpha \subseteq G$ denote the corresponding {\it root subgroup}, see \cite[(II.1.2)]{Ja3}. 

\bigskip

\begin{Corollary} \label{VV2}  Let $\lambda \in X(T)$ be a weight such that $\Psi_\lambda \ne \Psi$. Given $r \le s$, we have $\Sigma \cap (\Psi\setminus \Psi_\lambda) \neq \emptyset$ 
as well as 
\[ \cV_{({U_\alpha})_r}(k) \subseteq \cV_{G_r}(Z_s(\lambda)) \ \ \ \ \ \ \forall \ \alpha \in \Sigma \cap (\Psi\setminus \Psi_\lambda).\]
In particular, $\cx_{G_r}(Z_s(\lambda)) \ge r$. \end{Corollary}

\begin{proof} Since $\Psi_\lambda \subseteq \Psi$ is a subsystem (\cite[(2.7)]{Ja2}), there exists a simple root $\alpha$ in $\Psi \setminus \Psi_\lambda$. Given such a root $\alpha$, we 
let $L\supseteq T$ be the Levi subgroup of $G$ defined by $\alpha$. Thanks to Proposition \ref{Pr2}, the variety $\cV_{G_s}(Z_s(\lambda))$ contains $\cV_{L_s}(Z_s^L(\lambda))$, so that
\[ \cV_{L_r}(Z^L_s(\lambda)) \subseteq \cV_{G_r}(Z_s(\lambda)).\]
By assumption, the weight $\lambda$ is $p$-regular on $L$. Observe that the root group $U_\alpha$ is the unipotent radical of the Borel subgroup $B \cap L$ of $L$. Moreover, $U_\alpha$  is 
isomorphic to the additive group $\GG_a \cong {\rm Spec}_k(k[T])$. Accordingly, we have $(U_\alpha)_r \cong \GG_{a(r)}$, whence 
\[\cV_{(U_\alpha)_r}(k) \cong \cV_{\GG_{a(r)}}(k) \cong k^{r}\]
is irreducible, see \cite[(4.1)]{CPSvdK}. Consequently, Lemma \ref{VV1} yields $\cV_{L_r}(Z^L_s(\lambda)) = \cV_{(U_\alpha)_r}(k)$, whence $\cV_{(U_\alpha)_r}(k) \subseteq 
\cV_{G_r}(Z_s(\lambda))$ and $\cx_{G_r}(Z_s(\lambda)) \ge \dim \cV_{(U_\alpha)_r}(k) = r$. \end{proof}

\bigskip

\begin{Corollary} \label{VV3} Suppose that $p\ge 5$. Let $\lambda \in X(T)$ be a weight such that 

{\rm (a)} \ $\Psi_\lambda \ne \Psi$, and

{\rm (b)} \ $\cx_{G_r}(Z_r(\lambda)) = r$.

\noindent
Then there exists a simple root $\alpha \in \Sigma$ such that $\cV_{G_r}(Z_r(\lambda)) = \cV_{(U_\alpha)_r}(k)$ is irreducible.\end{Corollary}

\begin{proof} 
According to Corollary \ref{VV2}, we can find a simple root $\alpha \in \Psi$ such that the $r$-dimensional irreducible variety $\cV_{(U_\alpha)_r}(k)$ is contained in 
$\cV_{G_r}(Z_r(\lambda))$. By the arguments of Lemma  \ref{AR-V1}, this implies
\[ V_r(U_\alpha) \subseteq V_r(G)_{Z_r(\lambda)}.\]
In view of (b), the rank variety $V_r(U_\alpha)$ is an irreducible component of the $B$-invariant variety $V_r(G)_{Z_r(\lambda)}$. Thus, the stabilizer $H \subseteq B$ of $V_r(U_\alpha)$ is 
a closed subgroup of finite index, and the connectedness of $B$ gives $H = B$, so that the variety $V_r(U_\alpha)$ is also $B$-invariant. In particular,  the $B$-orbit of the canonical 
isomorphism $\GG_{a(r)} \stackrel{\sim}{\lra} (U_\alpha)_r$ is contained in $V_r(U_\alpha)$, implying the $B$-invariance of $(U_\alpha)_r \subseteq B_r$. 
Thus, $\fg_\alpha = \Lie((U_\alpha)_r)$ enjoys the same property, and the proof of \cite[(3.3)]{FR} provides a decomposition $G = SH$ of $G$ as an almost direct product with

(1) \ $\fg_\alpha \subseteq \Lie (S)$, and

(2) \ $S$ is almost simple of  type $A_1$, and

(3) \ $\fg = \Lie(S)\oplus \Lie(H)$.

\noindent
In view of (3), we have $(S\cap H)_1 \subseteq {\rm Cent}(S)_1 = e_k$, so that $(S\cap H)_r = e_k$.  Thus, the canonical map $S\times H \lra G$ induces a closed embedding $S_r \times H_r \hookrightarrow G_r$. Applying \cite[(I.9.6)]{Ja3} several times, we obtain
\[ \dim_k k[G_r] = p^{r\dim G} = p^{r(\dim S + \dim H)} = \dim_k k[S_r] \dim_k k[H_r] = \dim_k k[S_r\times H_r],\]
so that the above map is in fact an isomorphism. 

In view of \cite[(I.7.9(3))]{Ja3}, the arguments of \cite[(1.1)]{FR} show that the decomposition $G_r = S_r \times H_r$ induces an isomorphism 
\[ Z_r(\lambda) \cong Z_r^S(\lambda')\!\otimes_k\!Z_r^H(\lambda'')\]
between $Z_r(\lambda)$ and the outer tensor product of the corresponding Verma modules of the factors. The formula for varieties of outer tensor products, established in the proof of  
\cite[(7.2)]{SFB2}, implies
\[ (\ast) \ \ \ \ \ \ \ \ \ \ \cV_{G_r}(Z_r(\lambda)) \cong \cV_{S_r}(Z_r^S(\lambda')) \times \cV_{H_r}(Z_r^H(\lambda'')).\]
Since $\fg_\alpha \subseteq \Lie (S)$, we see that $\alpha$ is a root of $S$. Consequently, $U_\alpha \subseteq S$, and Corollary \ref{VV2} gives
\[\cV_{(U_\alpha)_r}(k)  \subseteq \cV_{S_r}(Z_r^S(\lambda')) \subseteq \cV_{(U_\alpha)_r}(k).\] 
Consequently, $\dim \cV_{S_r}(Z_r^H(\lambda)) = 0$, so that the conical variety $\cV_{S_r}(Z_r^H(\lambda))$ is a singleton. The isomorphism ($\ast$) now implies that the variety 
$\cV_{G_r}(Z_r(\lambda))$ is irreducible. Hence the inclusion $\cV_{(U_\alpha)_r}(k) \subseteq \cV_{G_r}(Z_r(\lambda))$ is in fact an equality. \end{proof} 

\bigskip

\begin{Remarks} (1) The case $r=1$ is of course a direct consequence of Carlson's Theorem \cite[Thm.~1]{Ca2}, which also holds in our context (cf.\ \cite[(7.7)]{SFB2}).

(2) If $G$ is almost simple, then the proof of (\ref{VV3}) yields that $G$ is of type $A_1$. In that case, Corollary \ref{VV2} and Proposition \ref{Pr2} show that $Z_r(\lambda)$ has an 
irreducible support variety of dimension $r$, whenever $\Psi_\lambda \ne \Psi$. 

(3) The reader is referred to the example at the end of Section \ref{s:depth} for a discussion of periodic baby Verma modules. \end{Remarks}

\bigskip
\noindent
Given $\lambda \in X(T)$ with $\langle \lambda+\rho,\alpha^\vee\rangle \ne 0$ for some $\alpha \in \Psi$, we define the \emph{depth} of $\lambda$ via
\[ \depth(\lambda) := \min\{s \in \NN_0 \mid \Psi_\lambda^s \ne \Psi\},\]
and put $\depth(\lambda) = -\infty$ otherwise. Let $\lambda \in X(T)$ and denote by $\cB_r(\lambda)$ the set of weights of the simple $\Dist(G_r)$-modules belonging to the block 
containing the simple $\Dist(G_r)$-module $L_r(\lambda)$. Thanks to \cite[(5.5)]{Ja1} (see also \cite[(II.9.14), (II.9.22)]{Ja3}), we have
\[ (\dagger) \ \ \ \ \ \ \ \ \  \cB_r(\lambda) = W\dact \lambda + p^{\depth(\lambda)} \ZZ \Psi + p^r X(T).\]
(Here we employ the convention $p^{-\infty} = 0$.) The notion of the depth of a weight will also play an important role in Section \ref{s:depth} below.

\bigskip

\begin{Proposition} \label{VV4} Suppose that $G$ is defined over $\FF_p$. Let $r \le s$ and set $t = {(s-r)\dim \mathfrak u}$. For every  weight $\lambda = \lambda_0 + p^r\lambda_1 
\in X(T)$ with $\lambda_0 \in X_r(T)$, the $G_rT$-module $\widehat{Z}_s(\lambda) |_{G_rT}$ possesses a filtration 
\[ (0) = M_0 \subseteq M_1 \subseteq \cdots \subseteq M_{p^t} = \widehat{Z}_s(\lambda) |_{G_rT} \]
by $G_rT$-submodules, such that $M_i/M_{i-1} \cong \widehat{Z}_r(\gamma_i)$ for some $\gamma_i \in \cB_r(\lambda_0)$. \end{Proposition}

\begin{proof} We first assume $(G,G)$ to be simply connected, so that every weight $\gamma \in X(T)$ is of the form
\[ \gamma = \gamma_0 + p^r\gamma_1,\]
with $\gamma_0 \in X_r(T)$ and $\gamma_1 \in X(T)$ (cf.\ \cite[(II.9.14)]{Ja3}). 

Since $\widehat{Z}_s(\lambda)$ is $B_r^-$-projective, \cite[(II.11.2)]{Ja3} provides a $\widehat{Z}_r$-filtration of $\widehat{Z}_s(\lambda) |_{G_rT}$ of length $p^t$. Let $\gamma \in 
X(T)$ be a weight, written as 
\[ \gamma = \mu + p^s\nu\]
with $\mu \in X_s(T)$. In view of \cite[(II.9.6(6))]{Ja3}, we have
\[ \widehat{L}_s(\gamma) \cong \widehat{L}_s(\mu)\!\otimes_k\!k_{p^s\nu}.\]
We now decompose $\mu \in X_s(T)$ as $\mu = \gamma_0 + p^r\gamma_1$ with $\gamma_0 \in X_r(T)$ and $\gamma_1 \in X(T)$. The assumption $\langle 
\gamma_1,\alpha^\vee\rangle < 0$ for some $\alpha \in \Sigma$ implies $\langle \mu, \alpha^\vee \rangle = \langle \gamma_0,\alpha^\vee \rangle +  p^r \langle 
\gamma_1,\alpha^\vee\rangle < \langle \gamma_0,\alpha^\vee \rangle - p^r <0$, a contradiction. Consequently, $\gamma_1 \in X(T)_+$, and Steinberg's tensor product theorem 
\cite[(II.3.16)]{Ja3} in conjunction with \cite[(II.9.6g)]{Ja3} yields 
\[ \widehat{L}_s(\mu) \cong L(\mu)|_{G_sT} \cong L(\gamma_0)|_{G_sT}\! \otimes_k\! L(\gamma_1)^{[r]}|_{G_sT}.\]
Upon restriction to $G_rT$, we obtain, observing \cite[(II.9.6(6))]{Ja3},
\begin{eqnarray*}
\widehat{L}_s(\gamma)|_{G_rT} & \cong  & \widehat{L}_s(\mu)|_{G_rT}\!\otimes_k\!k_{p^s\nu}|_{G_rT} \cong L(\gamma_0)|_{G_rT}\!\otimes_k\! L(\gamma_1)^{[r]}|_{G_rT}\!
\otimes_k\! k_{p^s\nu}|_{G_rT} \\
& \cong & \bigoplus_{\zeta \in X(T)} n_\zeta \widehat{L}_r(\gamma_0)\! \otimes_k\! k_{p^r\zeta+p^s\nu} \cong
\bigoplus_{\zeta \in X(T)} n_\zeta \widehat{L}_r(\gamma_0+p^r(\zeta+p^{s-r}\nu))\\
& \cong & \bigoplus_{\zeta \in X(T)} m_\zeta \widehat{L}_r(\gamma_0+p^r\zeta),
\end{eqnarray*}
for some $n_\zeta, m_\zeta \in \NN_0$. 

Now let $\eta \in X_r(T)$ be an arbitrary weight. Thanks to \cite[(II.11.2)]{Ja3}, the factor $\widehat{Z}_r(\eta)$ occurs in the $\widehat{Z}_r$-filtration of $\widehat{Z}_s(\lambda) 
|_{G_rT}$ with multiplicity
\[ \ell_\eta := \dim_k \Hom_{G_rT}(\widehat{Z}_s(\lambda)|_{G_rT} , \widehat{Z}_r'(\eta)).\]
If $\widehat{L}_s(\gamma)$ is a composition factor of the $G_sT$-module $\widehat{Z}_s(\lambda)$, then $L_s(\gamma)$ is a composition factor of the indecomposable $G_s$-module 
$Z_s(\lambda)$, so that $\gamma$ belongs to $\cB_s(\lambda)$. In view of the above isomorphisms, left exactness of $\Hom_{G_rT}( - , \widehat{Z}_r'(\eta))$ thus yields the 
estimate
\[ \ell_\eta \le \sum_{\gamma \in \cB_s(\lambda)} \sum_{\zeta \in X(T)} m_{\zeta,\gamma} \dim_k \Hom_{G_rT}(\widehat{L}_r(\gamma_0 + p^r\zeta), \widehat{Z}_r'(\eta)).\]
According to \cite[(II.9.6(2))]{Ja3}, $\ell_\eta \ne 0$ implies that $\eta = \gamma_0 + p^r\zeta$ for some $\gamma \in \cB_s(\lambda), \zeta \in X(T)$.

By Jantzen's description ($\dagger$) of $\cB_s(\lambda)$ above, there exist $w \in W$, $\alpha_0 \in \ZZ \Psi$, and $\chi \in X(T)$ such that
\[ \gamma = w\dact \lambda + p^{\depth(\lambda)}\alpha_0 + p^s\chi =  w\dact \lambda_0 + p^r w(\lambda_1) + p^{\depth(\lambda)}\alpha_0 + p^s\chi.\] 
Since $\lambda_0 \in X_r(T)$, we have $\depth(\lambda_0) \le r+1$. Consequently,
\[ \langle\lambda + \rho, \beta^\vee \rangle = \langle\lambda_0 + \rho, \beta^\vee \rangle + p^r \langle \lambda_1,\beta^\vee \rangle \in p^m\ZZ + p^r\ZZ \subseteq p^m \ZZ\]
for all $m \le \depth(\lambda_0)-1$ and $\beta \in \Psi$. Thus, we have $\depth(\lambda_0) \le \depth(\lambda)$ or $\depth(\lambda) = -\infty$.

Accordingly, setting $\beta_0: = p^{\depth(\lambda) - \depth(\lambda_0)}\alpha_0$ (with $\beta_0 := 0$ for $\depth(\lambda) = -\infty$), the above identity attains the form 
\[\gamma = w\dact \lambda_0 + p^{\depth(\lambda_0)}\beta_0 + p^r \omega,\]
for some $\omega \in X(T)$. This, however, implies that $\gamma = w\dact \lambda_0 + p^{\depth(\lambda_0)}\beta_0 + p^r \omega$ belongs to $\cB_r(\lambda_0)$. As a result, 
$\eta = \gamma_0 + p^r\zeta = \gamma + p^r(\zeta-\gamma_1)$ enjoys the same property. 

In general, there exists a covering $\tilde{G}$ of $G$ such that the derived group $(\tilde{G},\tilde{G})$ is simply connected. Our result now follows from Lemma \ref{AR-GT5}. \end{proof}

\bigskip

\begin{Corollary} \label{VV5} Suppose that $G$ is defined over $\FF_p$. Let $r \le s$ and $\lambda = \lambda_0 + p^r\lambda_1$ be a weight with $\lambda_0 \in X_r(T)$. The module 
$Z_s(\lambda)|_{G_r}$ belongs to the block $\cB_r(\lambda_0)$, and $\cV_{G_r}(Z_s(\lambda)) \subseteq \bigcup_{\gamma \in \cB_r(\lambda_0)} \cV_{G_r}(Z_r(\gamma))$.
\end{Corollary}

\begin{proof} This follows from Proposition \ref{VV4} and standard properties of support varieties, see \cite[(2.2.7)]{NPV}. \end{proof}

\bigskip
\noindent
In the following, we let 
\[ \St_r := \Dist(G_r)\!\otimes_{\Dist(B_r)}\!k_{(p^r -1)\rho}\]
denote the \emph{Steinberg module} of $G_r$. Our next result generalizes \cite[(II.11.8)]{Ja3}.

\bigskip

\begin{Corollary} \label{VV6} Suppose that $G$ is defined over $\FF_p$. Let $r \le s$. Then the following statements are equivalent:

{\rm (1)} \ $\Psi_\lambda^r = \Psi$.

{\rm (2)} \ There exists a simple projective $G_r$-module $P$ such that $Z_s(\lambda)|_{G_r} \cong p^{(s-r)\dim_k\fu} \, P$.

{\rm (3)} \  $Z_s(\lambda)|_{G_r}$ is projective. \end{Corollary}

\begin{proof} In view of Lemma \ref{AR-GT5}, we may assume that $(G,G)$ is simply connected.

(1) $\Rightarrow$ (2). Suppose that $\Psi_\lambda^r = \Psi$. Writing $\lambda = \lambda_0 + p^r\lambda_1$, with $\lambda_0 \in X_r(T)$, we obtain $\langle
\lambda_0+\rho,\alpha^\vee\rangle = p^r$ for all $ \alpha \in \Sigma$. Thanks to \cite[(II.11.8)]{Ja3}, this implies that $Z_r(\lambda_0)$ is a simple projective module and 
$\cB_r(\lambda_0) = \lambda_0 + p^rX(T)$. Thus, Proposition \ref{VV4} in conjunction with \cite[(II.3.7.(9))]{Ja3} provides a filtration of $Z_s(\lambda)|_{G_r}$ 
with $Z_r(\lambda_0)$ being the only filtration factor. Consequently,  $Z_s(\lambda)|_{G_r} \cong p^{(s-r)\dim_k\fu}\, Z_r(\lambda_0)$.

(2) $\Rightarrow$ (3). This is trivial.

(3) $\Rightarrow$ (1). Let $\alpha \in \Sigma$ be a simple root of $\Psi$ and let  $L$ be the corresponding Levi subgroup of $G$. By Proposition \ref{Pr2} we have
\[ \cV_{L_r}(Z_s^L(\lambda)) \subseteq \cV_{G_r}(Z_s(\lambda)).\]
Our current assumption thus implies the projectivity of the module $Z^L_s(\lambda)|_{L_r}$, see Proposition \ref{Pr1}. Thanks to Corollary \ref{VV5}, this module belongs to a block $\cB 
\subseteq \Dist(L_r) \cong \Dist(\SL(2)_r)^{r(\dim T-1)}$. If $\cB$ is not a simple block, then \cite[(II.11.16(b))]{Ja3} and \cite[(II.11.10),(II.11.11)]{Ja3} in conjunction with 
$\SL(2)_1$-theory imply that the dimensions of the principal indecomposable $\cB$-modules are even. As $\dim_kZ_s^L(\lambda) = p^s$, we have reached a contradiction. It follows that 
$Z^L_s(\lambda)|_{L_r} \cong p^{s-r} (\St_r\!\otimes_k\! k_\mu)$, where $\mu \in X(T)$ satisfies $\langle \mu,\alpha^\vee\rangle = 0$. Thus, $\langle \lambda + \rho, \alpha^\vee 
\rangle \in p^r\ZZ$, so that $\alpha \in \Psi^r_\lambda$. Since $\Psi_\lambda^r$ is a subsystem of $\Psi$ (\cite[(2.7)]{Ja2}), we conclude that $\Psi = \Psi_\lambda^r$, as desired. 
\end{proof}

\bigskip

\section{Depth Reduction for Verma Modules} \label{s:depth}
Throughout this section, $G$ is assumed to be a connected reductive smooth algebraic group that is defined over the Galois field $\FF_p$. 

We fix $r \in \NN$ and study non-projective Verma modules. Accordingly, we assume that $\lambda \in X(T)$ satisfies  $\Psi \ne \Psi^r_\lambda$, that is, $\depth(\lambda) \le r$. Recall 
that for any $d \in \NN$ the Steinberg module $\St_d \cong L((p^d-1)\rho)$ has the structure of a $G$-module, so that it is also a $G_rT$-module, see \cite[(II.10.1)]{Ja3}. Following 
\cite[(II.3.16)]{Ja3}, we denote by $M^{[d]}$ the \emph{$d$-Frobenius twist} of a module $M$, defined over a subgroup scheme $H \subseteq G$. As before, the characteristic $p$ is 
assumed to be good for $G$.

\bigskip

\begin{Proposition} \label{DR1} Let $1 \le d \le r$. The Frobenius homomorphism $F^d : G_r \lra G_{r-d}$ induces an isomorphism $\cV_{G_r}(\St_d) \stackrel{\sim}{\lra} 
\cV_{G_{r-d}}(k)$ that maps $\cV_{G_r}(M^{[d]}\!\otimes_k\!\St_d)$ onto $\cV_{G_{r-d}}(M)$ for every $M \in \modd G_{r-d}$. \end{Proposition}

\begin{proof} Given a $\Dist(G_r)$-module algebra $\Lambda$, we let $\HH^\bullet(G_r,\Lambda) := \bigoplus_{n\ge 0}\HH^{2n}(G_r,\Lambda)$ be the corresponding even 
cohomology ring. The algebras we are interested in are of the form
\begin{eqnarray*} \Lambda_d(M) & := &\Hom_k(M^{[d]}\!\otimes_k\!\St_d,M^{[d]}\!\otimes_k\!\St_d) \cong \Hom_k(M^{[d]},M^{[d]})\!\otimes_k\!\Hom_k(\St_d,\St_d)\\
 &\cong& \Hom_k(M,M)^{[d]}\!\otimes_k\!\Hom_k(\St_d,\St_d), \end{eqnarray*}
where $M \in \modd G_{r-d}$. Since $\Lambda_d(M)$ is a projective $G_d$-module, the Lyndon--Hochschild--Serre spectral sequence $\HH^n(G_r/G_d,\HH^m(G_d, \Lambda_d(M))) 
\Rightarrow \HH^{n+m}(G_r,\Lambda_d(M))$, cf.\ \cite[(I.6.6)]{Ja3}, collapses to isomorphisms
\[ \HH^n(G_r/G_d,\Lambda_d(M)^{G_d}) \cong \HH^n(G_r,\Lambda_d(M)) \ \ \ \ \forall \ n \ge 0.\]
In view of $\St_d = L((p^d-1)\rho)$ being a simple $G$-module, the groups $G_r$ and $G_r/G_d$ operate trivially on the one-dimensional space $\Hom_{G_d}(\St_d,\St_d)$. As
$G_d$ acts trivially on $\Hom_k(M,M)^{[d]}$, we obtain $G_r$-isomorphisms
\[ \Lambda_d(M)^{G_d} \cong \Hom_k(M,M)^{[d]}\!\otimes_k\! \Hom_{G_d}(\St_d,\St_d) \cong \Hom_k(M,M)^{[d]}.\]
The canonical quotient map $\pi : G_r \lra G_r/G_d$ and the natural map $\Hom_k(M,M)^{[d]} \lra \Lambda_d(M) \ ; $ $f \mapsto f \otimes \id_{\St_d}$ induce homomorphisms 
\[ \pi^\bullet_M : \HH^\bullet(G_r/G_d,\Hom_k(M,M)^{[d]}) \lra \HH^\bullet(G_r,\Hom_k(M,M)^{[d]}) \] 
and 
\[\Upsilon_M : \HH^\bullet(G_r,\Hom_k(M,M)^{[d]}) \lra \HH^\bullet(G_r,\Lambda_d(M))\] 
of $k$-algebras, respectively. By the above, the inflation map $\Upsilon_M \circ \pi^\bullet_M$ is bijective (cf.\ \cite[(6.7.3), (6.8.2)]{Wei}). Thus, $\pi^\bullet_M$ is injective and, 
specializing $M=k$, we obtain $\HH^\bullet(G_r,k) = \im \pi^\bullet_k \oplus \ker \Upsilon_k$. Consequently,
\[ \bar{\pi}^\bullet_k : \HH^\bullet(G_r/G_d,k) \lra \HH^\bullet(G_r,k)/\ker \Upsilon_k \ \ ; \ \  x \mapsto \pi^\bullet_k(x) + \ker \Upsilon_k\]
is an isomorphism of $k$-algebras. Thanks to \cite[(I.9.5)]{Ja3}, the Frobenius endomorphism $F^d$ induces an isomorphism $\omega_d : G_r/G_d \lra G_{r-d}$, so that 
\[ \bar{\pi}^\bullet_k \circ \omega_d^\bullet : \HH^\bullet(G_{r-d},k) \lra \HH^\bullet(G_r,k)/\ker\Upsilon_k\]
is also an isomorphism. By definition, the variety $\cV_{G_r}(\St_d)$ is the zero locus of the ideal $\ker\Upsilon_k$, so that the associated comorphism is the desired isomorphism
\[ \cV_{G_r}(\St_d) \stackrel{\sim}{\lra} \cV_{G_{r-d}}(k).\]
(By construction, this isomorphism is the restriction of the morphism $\cV_{G_r}(k) \lra \cV_{G_{r-d}}(k)$ induced by $F^d$.) The Frobenius homomorphism $F^d$ gives rise to a 
commutative diagram
\[ \begin{CD} \HH^\bullet(G_{r-d},k) @>\Phi_M>> \HH^\bullet(G_{r-d},\Hom_k(M,M)) \\
        @V\pi^\bullet_k \circ \omega_d^\bullet VV @V\pi^\bullet_M\circ \omega_dVV\\
 \HH^\bullet(G_r,k) @>\Phi_{M^{[d]}}>> \HH^\bullet(G_r,\Hom_k(M,M)^{[d]}), \end{CD} \]
where the horizontal maps are given by $\alpha \mapsto \alpha.\id_M$. Since $\Upsilon_M \circ \pi_M^\bullet$ is bijective, we conclude the injectivity of the right-hand map. Owing 
to \cite[(7.2)]{SFB2}, we have
\[\cV_{G_r}(M^{[d]}\!\otimes_k\!\St_d) = \cV_{G_r}(M^{[d]})\cap \cV_{G_r}(\St_d) \subseteq \cV_{G_r}(\St_d).\]
The intersection is given by the maximal ideal spectrum of $\HH^\bullet(G_r,k)/I$, where $I = \ker \Phi_{M^{[d]}}+\ker \Upsilon_k$. From the above diagram, we conclude that the image 
of $I$ in $\HH^\bullet(G_r,k)/\ker \Upsilon_k$ corresponds under the isomorphism $\bar{\pi}^\bullet_k\circ \omega_d^\bullet$ to the ideal $\ker\Phi_M$. Consequently, the 
isomorphism $\cV_{G_r}(\St_d) \stackrel{\sim}{\lra} \cV_{G_{r-d}}(k)$ induced by $\bar{\pi}^\bullet_k \circ \omega_d^\bullet$, sends $\cV_{G_r}(M^{[d]})\cap \cV_{G_r}(\St_d)$ 
onto $\cV_{G_{r-d}}(M)$, as desired. \end{proof}

\bigskip

\begin{Theorem} \label{DR2} Suppose $G$ is semi-simple and simply connected and that $\lambda \in X(T)$ satisfies $r \ge \depth(\lambda) = d+1 \ge 2$. Then there exists a weight 
$\mu \in X(T)$ of depth $\depth(\mu) = 1$ such that

{\rm (a)} \ there is an isomorphism  $\widehat{Z}_r(\lambda) \cong \widehat{Z}_{r-d}(\mu)^{[d]}\!\otimes_k\!\St_d$, and

{\rm (b)} \ there is an isomorphism $\cV_{G_r}(Z_r(\lambda)) \cong \cV_{G_{r-d}}(Z_{r-d}(\mu))$. \end{Theorem}

\begin{proof} Let $e_d \in \Dist(G_d) \subseteq \Dist(G_r)$ be the central primitive  idempotent, defined by the projective simple $G_d$-module $\St_d$. Owing to \cite[(II.10.3)]{Ja3}, 
$e_d$ is a central idempotent of $\Dist(G_r)$. Since $\depth(\lambda) = d+1$, we have $\Psi_\lambda^d = \Psi$. As $G$ is semi-simple and simply connected, $\St_d$ is the only simple, 
projective $G_r$-module. Accordingly, Corollary \ref{VV6} provides an isomorphism  $Z_r(\lambda)|_{G_d} \cong n \St_d$, so that $e_dZ_r(\lambda) = Z_r(\lambda)$. Thanks to \cite[(II.10.4.b)]{Ja3}, there exists a $G_r T/G_d$-module $M$ such that 
\[ Z_r(\lambda) \cong M\!\otimes_k\!\St_d.\]
In particular, $\dim_k M = p^{(r-d)\dim_k \fu}$.

Setting $H = G_d B_r^-$, we propose to show that $Z_r(\lambda)$ is a projective $H$-module. Recall that $V_r(H)_{Z_r(\lambda)}  \subseteq V_r(G)_{Z_r(\lambda)} \subseteq V_r(U)$, 
whence
\[ V_r(H)_{Z_r(\lambda)} = V_r(H\cap U_r)_{Z_r(\lambda)}.\]
Let $R$ be a commutative $k$-algebra and let $x$ be an element of $(H\cap U_r)(R)$. By general theory \cite[(15.5)]{Wa}, there exists a faithfully flat extension $S$ of $R$ such that $x_S = 
ab$ for $a \in G_d(S)$ and $b \in B_r^-(S)$. Here $x_S$ denotes the image of $x$ in $(H\cap U_r)(S)$ under the canonical map induced by $R \lra S$. An application of the $d$-th power of 
the Frobenius endomorphism gives
\[ F^d(x)_S = F^d(x_S) = F^d(a)F^d(b) = F^d(b) \in U_{r-d}(S)\cap B^-_{r-d}(S) = \{1\}.\]
In view of \cite[(15.6)]{Wa}, we also have $F^d(x) = 1$, whence $x \in U_d(R)$. Thus, $H\cap U_r \subseteq U_d$ implying
\[V_r(H\cap U_r)_{Z_r(\lambda)} \subseteq V_r(U_d)_{Z_r(\lambda)} \subseteq V_r(G_d)_{Z_r(\lambda)} \cong V_d(G_d)_{Z_r(\lambda)},\]
where the last isomorphism is obtained by restricting the canonical identification $V_d(G_d) \lra V_r(G_d) \ \ ; \ \varphi \mapsto \varphi \circ F^{r-d}$.  In view of Corollary \ref{VV6}, the 
right-hand variety is trivial, proving that $V_r(H)_{Z_r(\lambda)} = V_r(H\cap U_r)_{Z_r(\lambda)} = \{0\}$. As a result, the module $Z_r(\lambda)|_H$ is projective. 

Owing to \cite[(I.9.5)]{Ja3}, the Frobenius homomorphism $F^d$ induces an isomorphism $G_rT/G_d \cong G_{r-d}T$. Consequently, $M = N^{[d]}$ for some $G_{r-d}T$-module $N$. 
Moreover, $N$ is a projective $B_{r-d}^-$-module if and only if $M$ is a projective $H/G_d$-module. By \cite[(II.10.4)]{Ja3}, the functor $X \mapsto X\!\otimes_k\!\St_d$ gives rise to an 
equivalence between $\modd \Dist(H/G_d)$ and  $\modd \Dist(H) e_d$, sending $M$ to $Z_r(\lambda)$. Since $Z_r(\lambda)$ is a projective $H$-module, we conclude the 
$B_{r-d}^-$-projectivity of $N$. According to \cite[(II.11.2)]{Ja3}, the $G_{r-d} T$-module $N$ therefore affords a filtration by baby Verma modules. Thus, for dimension reasons, we have 
$N \cong Z_{r-d}(\nu)$ for some $\nu \in X(T)$.

Since
\[ \fF(\widehat{Z}_{r-d}(\nu)^{[d]}\!\otimes_k\!\St_d) \cong \fF(\widehat{Z}_{r-d}(\nu)^{[d]})\!\otimes_k\!\fF(\St_d) \cong Z_{r-d}(\nu)^{[d]}\!\otimes_k\!\St_d \cong 
Z_r(\lambda),\]
it follows from \cite[(4.1)]{GG1} that there exists $\gamma \in X(T)$ such that
\[ \widehat{Z}_r(\lambda) \cong \widehat{Z}_{r-d}(\nu)^{[d]}\!\otimes_k\!\St_d \otimes_k\!k_{p^r\gamma} \cong 
\widehat{Z}_{r-d}(\nu+p^{r-d}\gamma)^{[d]}\!\otimes_k\!\St_d.\]
A comparison of highest weight vectors gives $\lambda = p^d \nu +p^r\gamma + (p^d -1)\rho$. Consequently, 
\[\langle \lambda + \rho, \alpha^\vee \rangle = p^d \langle \nu + p^{r-d}\gamma+ \rho, \alpha^\vee \rangle\] 
for all $\alpha \in \Psi$, so that the weight $\mu := \nu+p^{r-d}\gamma$ has depth $\depth(\mu) = 1$. This completes the proof of (a). Part (b) is now a direct consequence of
Proposition \ref{DR1}. \end{proof}

\bigskip

\begin{Example} We consider the semi-simple, simply connected group $G = \SL(2)$. If $\lambda \in X(T)$ is a weight with $\depth(\lambda) \le r$, then
\[ \cV_{\SL(2)_r}(Z_r(\lambda)) \cong k^{r+1-\depth(\lambda)}.\]
If $\depth(\lambda) = 1$, then our assertion follows directly from Corollary \ref{VV2}. Alternatively, Theorem \ref{DR2} reduces us to this case. \end{Example}

\bigskip

\begin{Corollary} \label{DR3} Let $p\ge h$ and suppose that $G$ is semi-simple and simply connected. Let $\lambda \in X(T)$ be a weight of depth $\depth(\lambda) \le r$ such that 
$\cV_{U_{r+1-\depth(\lambda)}}(k)$ is irreducible. Then $\cV_{G_r}(Z_r(\lambda)) \cong \cV_{U_{r+1-\depth(\lambda)}}(k)$ if and only if $\lambda$ is 
$p^{\depth(\lambda)}$-regular. \end{Corollary}

\begin{proof} Thanks to Theorem \ref{DR2}, there exists a weight $\mu \in X(T)$ of depth $1$ and isomorphisms
\[ Z_r(\lambda) \cong Z_{r-d}(\mu)^{[d]}\!\otimes_k\!\St_d \ \ \text{and} \ \ \cV_{G_r}(Z_r(\lambda)) \cong \cV_{G_{r-d}}(Z_{r-d}(\mu))\]
of modules and varieties, respectively. According to Lemma \ref{VV1}, the condition $\cV_{G_{r-d}}(Z_{r-d}(\mu)) \cong \cV_{U_{r-d}}(k)$ is equivalent to $\mu$ being $p$-regular. Since 
$\lambda = p^d\mu +(p^d-1)\rho$, this happens precisely when $\lambda$ is $p^{\depth(\lambda)}$-regular. \end{proof}

\bigskip

\begin{Example} Consider the semi-simple, simply connected group $G = \SL(3)$. If $\lambda \in X(T)$ has depth $\depth(\lambda) \le r$ and $\lambda$ is 
$p^{\depth(\lambda)}$-regular, then the example following Lemma \ref{VV1} in conjunction with Corollary \ref{DR3} implies $\cV_{\SL(3)_r}(Z_r(\lambda)) \cong 
\cV_{U_{r+1-\depth(\lambda)}}(k)$ as well as $\cx_{\SL(3)_r}(Z_r(\lambda)) = 2(r-\depth(\lambda))+3$. \end{Example}

\bigskip

\begin{Corollary} \label{DR4} Let $p\ge h$ and suppose that $G$ is semi-simple and simply connected. Let $\lambda \in X(T)$ be a weight of depth $\depth(\lambda) = r$. Then 
$\cV_{G_r}(Z_r(\lambda)) \cong \cV_{U_1}(k)$ if and only if $\lambda$ is $p^r$-regular. \hfill $\square$ \end{Corollary}

\bigskip

\begin{Corollary} \label{DR5} Suppose that $G$ is semi-simple, simply connected and of characteristic $p\ge 5$. Let $\lambda \in X(T)$ be a weight such that $\cx_{G_r}(Z_r(\lambda)) = 
r+1-\depth(\lambda)$. Then there exists a simple root $\alpha \in \Psi$ such that $\cV_{G_r}(Z_r(\lambda)) \cong \cV_{(U_\alpha)_{r+1-\depth(\lambda)}}(k)$ is irreducible.
\end{Corollary}

\begin{proof} As before, we write $\depth(\lambda) = d+1$. If $\depth(\lambda) = 1$, then the result follows directly from Corollary \ref{VV3}. Alternatively, Theorem \ref{DR2}(b)
provides a weight $\mu$ of depth $1$ and an isomorphism $\cV_{G_r}(Z_r(\lambda)) \cong \cV_{G_{r-d}}(Z_{r-d}(\mu))$, so that another application of Corollary \ref{VV3} 
completes the proof.  \end{proof}

\bigskip
\noindent
The following result reduces questions concerning AR-components of Verma modules to those, whose defining weights have depth $1$. For ease of notation, we shall often denote Verma 
modules defined over different groups by the same symbol.

\bigskip

\begin{Theorem} \label{DR6} Let $\lambda \in X(T)$ be a weight such that $r \ge \depth(\lambda) = d+1 \ge 1$. Then there exists a covering group $\tilde{G}$ of $G$ and functors
 \[ \Phi : \modd G_rT \lra \modd \tilde{G}_{r-d}\tilde{T} \ \ \text{and} \ \  \Upsilon : \modd G_r \lra \modd \tilde{G}_{r-d}\]
 with the following properties:

 {\rm (a)} \ The diagram 
 \[ \begin{CD} \modd G_rT @> \Phi >> \modd \tilde{G}_{r-d}\tilde{T} \\
        @V\fF VV @V\tilde{\fF} VV\\
 \modd G_r @>\Upsilon >>  \mod \tilde{G}_{r-d} \end{CD} \]
 commutes.
 
 {\rm (b)} \ There is a weight $\mu \in X(\tilde{T})$ of depth $\depth(\mu) = 1$ such that $\Phi(\widehat{Z}_r(\lambda)) \cong \widehat{Z}_{r-d}(\mu)$.
 
 {\rm (c)} \ We have $\cx_{G_r}(\widehat{Z}_r(\lambda)) = \cx_{\tilde{G}_{r-d}}(\widehat{Z}_{r-d}(\mu)) \ge r-d$.
 
 {\rm (d)}  \ The functor $\Phi$ induces an isomorphism $\widehat{\Theta}_r(\lambda) \lra \widehat{\Theta}_{r-d}(\mu) \ \ ; \ \ [M] \mapsto [\Phi(M)]$ of stable translation
 quivers. \end{Theorem}

\begin{proof} In view of Corollary \ref{VV2}, our statement holds for $d=0$. We therefore assume that $d\ge 1$. 

Proceeding in several steps, we first prove our result under the assumption:

\medskip

(i) \ {\it $G$ is semi-simple and simply connected}. 

\smallskip
\noindent
Owing to Theorem \ref{DR2}, there exists a weight $\mu \in X(T)$ such that $\depth(\mu) = 1$ and
\[ \widehat{Z}_r(\lambda) \cong \widehat{Z}_{r-d}(\mu)^{[d]}\! \otimes_k\! \St_d.\]
Letting $e_d \in \Dist(G_d) \subseteq \Dist(G_r)$ be the primitive idempotent belonging to $\St_d$, we denote by $\modd G_rTe_d$ the full subcategory of $\modd G_rT$, consisting of 
those $G_rT$-modules $M$ satisfying $e_dM = M$. In view of \cite[(II.10.4)]{Ja3}, $\modd G_rTe_d$ is a sum of blocks of $\modd G_rT$.

Since the Frobenius endomorphism $F^d$ induces an isomorphism $G_rT/G_d \cong G_{r-d}T$, an application of \cite[(II.10.4/5)]{Ja3} shows that the restriction of the functor
\[ \Phi : \modd G_rT  \lra  \modd G_{r-d}T \ \ ; \ \  M \mapsto \Hom_{G_d}(\St_d,M)^{[-d]}\]
to $\modd G_rTe_d$ is an equivalence of categories. As such, $\Phi$ induces an isomorphism between the stable Auslander--Reiten quivers $\Gamma_s(\modd G_rTe_d)$ and 
$\Gamma_s(G_{r-d}T)$, sending the class of the module $\widehat{Z}_r(\lambda) \in \modd G_rTe_d$ onto that of $\widehat{Z}_{r-d}(\mu)$. Since $\modd G_rTe_d$ is a sum of 
blocks of $\modd G_rT$, the component of $\Gamma_s(\modd G_rTe_d)$ containing $[\widehat{Z}_r(\lambda)]$ coincides with $\widehat{\Theta}_r(\lambda)$. Consequently, 
$\Phi(\widehat{\Theta}_r(\lambda)) = \widehat{\Theta}_{r-d}(\mu)$, so that (b) and (d) hold for $\tilde{G} := G$.

Another application of \cite[(II.10.4/5)]{Ja3} shows that the restriction of
\[ \Upsilon : \modd G_r  \lra  \modd G_{r-d} \ \ ; \ \  M \mapsto \Hom_{G_d}(\St_d,M)^{[-d]},\]
to $\modd \Dist(G_r)e_d$ is an equivalence of categories. By definition, we have $\fF \circ \Phi = \Upsilon \circ \fF$. Since $\Psi_\mu \ne \Psi$, a consecutive application of Corollary 
\ref{VV2} and Theorem \ref{DR2} yields
\[ r-d \le \cx_{G_{r-d}}(\widehat{Z}_{r-d}(\mu)) = \cx_{G_r}(\widehat{Z}_r(\lambda)),\]
as desired. \hfill $\Diamond$

\medskip

(ii) \ {\it $G$ is semi-simple}.

\smallskip
\noindent
We consider a semi-simple simply connected covering group $\tilde{G}$ along with the maps $\pi: \tilde{G}_r \tilde{T} \lra G_rT$ and $\sigma : \tilde{G_r} \lra G_r$ that are induced by the 
canonical quotient map $\omega : \tilde{G} \lra G$. By Lemma \ref{AR-GT5}, the corresponding pull-back functor $\pi^\ast: \modd G_rT \lra \modd \tilde{G}_r\tilde{T}$ allows us to view 
$\modd G_rT$ as the sum of those blocks of $\modd \tilde{G}_r\tilde{T}$, whose modules have weights belonging to $X(T) \subseteq X(\tilde{T})$.

We first show that the pull-back $\pi^\ast(\widehat{Z}_r(\lambda))$ coincides with the baby Verma module $\widehat{\tilde{Z}}_r(\lambda)$ of $\modd \tilde{G}_r\tilde{T}$. Since 
$\tilde{G}$ is a covering, the map $\omega$ induces an isomorphism $\tilde{U}^- \lra U^-$, cf.\ \cite[(II.9.7)]{Ja3}. Consequently, $\pi^\ast(\widehat{Z}_r(\lambda))$ is a 
$\tilde{U}_r^-$-projective $\tilde{G}_r\tilde{T}$-module. In view of the isomorphism $\tilde{B}_r^- \cong \tilde{U}_r^-\rtimes \tilde{T}_r$, our module is even 
$\tilde{B}_r^-$-projective. Thus, \cite[(II.11.2)]{Ja3} provides a filtration of $\pi^\ast(\widehat{Z}_r(\lambda))$ by baby Verma modules, which, for dimension reasons, only has one 
constituent $\widehat{\tilde{Z}}_r(\mu)$. Since all weights of $\pi^\ast(\widehat{Z}_r(\lambda))$ are $\le \lambda$, while all weights of $\widehat{\tilde{Z}}_r(\mu)$ are $\le \mu$, 
we obtain $\lambda = \mu$.

The depth of $\lambda \in X(T)$ coincides with the depth of $\lambda$, viewed as an element of $X(\tilde{T})$ (see \cite[(II.1.17)]{Ja3}). By part (i), we have functors $\Phi_1 : \modd 
\tilde{G}_r\tilde{T}  \lra  \modd \tilde{G}_{r-d}\tilde{T}$ and $\Upsilon_1 : \modd \tilde{G}_r \lra \modd \tilde{G}_{r-d}$ satisfying (a)-(c). In particular, there exists a weight $\mu \in 
X(\tilde{T})$ of depth $1$ such that
\[ \Phi_1(\widehat{\tilde{Z}}_r(\lambda)) \cong \widehat{\tilde{Z}}_{r-d}(\mu).\]
We now consider the functors $\Phi_2 := \Phi_1\circ \pi^\ast$ and $\Upsilon_2 := \Upsilon_1\circ \sigma^\ast$. We thus have $\tilde{\fF} \circ \pi^\ast = \sigma^\ast \circ \fF$, 
whence
\[ \tilde{\fF} \circ \Phi_2 = \Upsilon_1\circ \tilde{\fF} \circ \pi^\ast = \Upsilon_1 \circ \sigma^\ast \circ \fF = \Upsilon_2\circ \fF.\]
Moreover, there are isomorphisms
\[ \Phi_2(\widehat{Z}_r(\lambda)) \cong \Phi_1(\pi^\ast(\widehat{Z}_r(\lambda))) \cong \Phi_1(\widehat{\tilde{Z}}_r(\lambda)) \cong \widehat{\tilde{Z}}_{r-d}(\mu),\]
so that (a) and (b) hold.

In view of (i), property (c) follows from the equality $\cx_{\tilde{G}_r}(\pi^\ast(\widehat{Z}_r(\lambda))) = \cx_{G_r}(\widehat{Z}_r(\lambda))$, which can be seen as follows: Let 
$(\widehat{P}_n)_{n\ge 0}$ be a minimal projective resolution of the $G_rT$-module $\widehat{Z}_r(\lambda)$. Owing to Lemma \ref{AR-GT5}, $(\pi^\ast(\widehat{P}_n))_{n\ge 0}$ is 
a minimal projective resolution of $\pi^\ast(\widehat{Z}_r(\lambda))$. Since $\Top_{G_r}(\fF(M)) = \fF(\Top_{G_rT}(M))$ for every $G_rT$-module $M$ (see \cite[(II.9.6.(11))]{Ja3}), 
the forgetful functor $\fF$ takes projective covers to projective covers (cf.\ also \cite[(1.3)]{GG2} and \cite[(1.6)]{Ja4}). Consequently, the complexities 
$\cx_{G_r}(\widehat{Z}_r(\lambda))$ and $\cx_{\tilde{G}_r}(\pi^\ast(\widehat{Z}_r(\lambda)))$ are computable from the growths of the sequences $(\dim_k \widehat{P}_n)_{n\ge 0} 
= (\dim_k \pi^\ast(\widehat{P}_n))_{n\ge 0}$. 

By virtue of Lemma \ref{AR-GT5}, the functor $\pi^\ast$ induces an isomorphism 
\[ \widehat{\Theta}_r(\lambda) \lra \widehat{\tilde{\Theta}}_r(\lambda) \ \ ; \ \ [M] \mapsto [\pi^\ast(M)],\]
so that $\Phi_2$ enjoys the required property (d). \hfill $\Diamond$

\medskip

(iii) \ {\it $G = G'\times T''$, with a semi-simple group $G'$ and a torus $T''$}.

\smallskip
\noindent
We denote by $T'$ a maximal torus of $G'$, so that $T := T'\times T''$ is a maximal torus of $G$. Then $X(T) \cong X(T')\times X(T'')$, and our weight $\lambda$ has the form
\[ \lambda = (\lambda',\lambda''),\]
where $\depth(\lambda') = \depth(\lambda)$. According to step (ii), we can find a covering $\tilde{G'} \lra G'$ and functors 
\[ \Phi_2 : \modd G'_rT' \lra \modd \tilde{G'}_{r-d}\tilde{T'} \ \ \text{and}\ \  \Upsilon_2 : \modd G'_r \lra \modd \tilde{G'}_{r-d}\]
as well as a weight $\mu' \in X(\tilde{T'})$ of depth $1$ such that (a)-(d) hold.
 
Let $\tilde{G} := \tilde{G'} \times T''$, so that $\tilde{G}$ is a covering of $G$. Since $\modd G_rT$ is a sum of copies of $\modd G'_rT'$, with each constituent being given by those $G_rT$-modules, whose weights belong to $X(T')\times \{\gamma\}$ for a fixed $\gamma \in X(T'')$, we consider the composite
\[ \Phi_3 := \iota_{\lambda''} \circ \Phi_2 \circ {\rm pr}_{\lambda''},\]
where
\[ \iota_{\lambda''} : \modd \tilde{G'}_{r-d}\tilde{T'} \lra {\modd \tilde{G}_{r-d}\tilde{T} \ \ \text{and} \ \ \rm pr}_{\lambda''} : \modd G_rT \lra \modd G'_rT'\]
are the inclusion functor and the projection functor induced by the decompositions of the relevant categories. Note that ${\rm pr}_{\lambda''}$ sends $\widehat{Z}_r(\lambda)$ to 
$\widehat{Z}_r(\lambda')$, so that
\[ \Phi_3(\widehat{Z}_r(\lambda)) \cong \widehat{Z}_{r-d}(\mu',\lambda'').\]
As noted earlier, the weight $\mu := (\mu',\lambda'')$ satisfies $\depth(\mu) = \depth(\mu') = 1$. Since the elements of $\widehat{\Theta}_r(\lambda)$ and 
$\widehat{\Theta}_{r-d}(\mu)$ belong to the respective blocks defined by $\lambda''$, the functors $\iota_{\lambda''}$ and ${\rm pr}_{\lambda''}$ give rise to isomorphisms 
$\widehat{\Theta}_{r-d}(\mu') \lra \widehat{\Theta}_{r-d}(\mu)$ and $\widehat{\Theta}_r(\lambda) \lra \Theta_r(\lambda')$ of stable translation quivers. Consequently, $\Phi_3$ 
satisfies (d). 

The module category $\modd G_r$ has a decomposition that is compatible with the one given above. Hence the definition
\[ \Upsilon_3 := \iota_{\lambda''} \circ \Upsilon_2 \circ {\rm pr}_{\lambda''}\]
yields a functor so that the pair $(\Phi_3,\Upsilon_3)$ satisfies (a).

 The block decompositions give rise to
 \[ \cx_{G_r}(\widehat{Z}_r(\lambda)) = \cx_{G'_r}(\widehat{Z}_r(\lambda')) \ \ \text{and} \ \ \cx_{\tilde{G}_{r-d}}(\widehat{Z}_{r-d}(\mu)) = 
 \cx_{\tilde{G'}_{r-d}}(\widehat{Z}_{r-d}(\mu')),\]
 so that (c) follows from (ii). \hfill $\Diamond$
 
 \medskip
 
 (iv) {\it $G$ is reductive}.
 
 \smallskip
 \noindent
General theory provides a covering $(G,G)\times Z(G)^\circ \stackrel{\omega}{\lra} G$, with the group $G' := (G,G)\times Z(G)^\circ$ being of the type discussed in (iii). Let $\tilde{G'}$ 
be a covering of $G'$ satisfying (a)-(d). Arguing as in step (ii), we obtain the desired functors by combining those of (iii) with the pull-back functors defined by $\omega$. Since the kernel of 
the canonical map $\tilde{G'} \lra G$ is a normal subgroup of dimension zero, it is contained in the center $Z(\tilde{G'})$ and thus in particular, it is abelian. As it also is an extension of two 
diagonalizable group schemes, \cite[(IV.\S1.(4.5))]{DG} ensures that it is diagonalizable. Consequently, the group $\tilde{G'}$ is indeed a covering of $G$. \end{proof}
 
\bigskip
\noindent
For our next result, recall the notation $x^+$ from Lemma \ref{AR-GT3} for the set of successors of a vertex $x$ in the stable AR-quiver $\Gamma_s(G_r)$.

\bigskip

\begin{Corollary} \label{DR7} Let $\lambda \in X(T)$ be a weight such that $r \ge \depth(\lambda) = d+1 \ge 1$.  Then there exists a covering group $\tilde{G}$ and a weight $\mu \in 
X(\tilde{T})$ with the following properties.
 
{\rm (a)} \ There exists a surjective morphism $\Upsilon : \Theta_r(\lambda) \lra \Theta_{r-d}(\mu)$ of stable translation quivers such that $\Upsilon(x^+) = \Upsilon(x)^+$ for all $x 
\in \Theta_r(\lambda)$.
 
 {\rm (b)} \ If $G$ is semi-simple and simply connected, then $\Upsilon$ is an isomorphism. \end{Corollary}
 
\begin{proof} According to Theorem \ref{DR6}, there exists a covering group $\tilde{G}$ and functors 
 \[ \Phi : \modd G_rT \lra \modd \tilde{G}_{r-d}\tilde{T} \ \ \text{and} \ \ \Upsilon : \modd G_r \lra \modd \tilde{G}_{r-d}\]
giving rise to a commutative diagram 
 \[ \begin{CD}  \modd G_rT @> \Phi >> \modd \tilde{G}_{r-d}\tilde{T} \\
        @V\fF VV @V\tilde{\fF} VV\\
 \modd G_r @>\Upsilon >>  \mod \tilde{G}_{r-d}. \end{CD} \]
 In addition, the functor $\Phi$ induces an isomorphism 
 \[\widehat{\Theta}_r(\lambda) \lra \widehat{\Theta}_{r-d}(\mu) \ \ ; \ \  [M] \mapsto [\Phi(M)]\]
 for a suitable weight $\mu \in X(\tilde{T})$ of depth $1$. In view of Lemma \ref{AR-GT3} and the above diagram, the map
 \[ \Upsilon : \Theta_r(\lambda) \lra \Theta_{r-d}(\mu) \ \ ; \ \ [M] \mapsto [\Upsilon(M)]\]
 is well-defined. 
 
If $\alpha : x \rightarrow y$ is an arrow in $\Theta_r(\lambda)$, then there exists an arrow $\hat{\alpha} : \hat{x} \rightarrow \hat{y}$ in $\widehat{\Theta}_r(\lambda)$ such that 
$\fF(\hat{x}) = x$ and $\fF(\hat{y}) = y$. Consequently, $\Upsilon(\alpha) := \tilde{\fF}(\Theta(\hat{\alpha}))$ is an arrow in $\Theta_{r-d}(\mu)$ with starting point 
$\tilde{\fF}(\Phi(\hat{x})) = \Upsilon(\fF(\hat{x})) = \Upsilon(x)$ and endpoint $\tilde{\fF}(\Phi(\hat{y})) = \Upsilon(y)$. With this definition, $\Upsilon$ becomes a morphism of 
quivers. Since $\fF$, $\Phi$, and $\tilde{\fF}$ commute with translations, so does $\Upsilon$. Finally, given $x = \fF(\hat{x}) \in \Theta_r(\lambda)$, we obtain, observing Lemma 
\ref{AR-GT3}, 
 \[\Upsilon(x^+) = \Upsilon(\fF(\hat{x}^+)) = \tilde{\fF}(\Phi(\hat{x}^+)) = \tilde{\fF}(\Phi(\hat{x})^+) = (\tilde{\fF}\circ \Phi)(\hat{x})^+ = \Upsilon(x)^+,\] 
 as desired.
  
 Property (b) follows directly from step (i) of the proof of Theorem \ref{DR6}. \end{proof}
 
 \bigskip
 
 \begin{Example} Let $G = \SL(2)$ and $\lambda = p^{r-1}a-1$, $a \not \in p\ZZ$ be a weight of depth $r$. If $r \ge 2$, then Corollary \ref{DR7} implies $\Theta_r(\lambda) \cong 
 \Theta_1(a-1)$. Since $\Theta_1(a-1)$ has a one-dimensional support variety, it is of the form $\ZZ[A_\infty]/(\tau)$, with $Z_1(a-1)$ having exactly one predecessor. Consequently, 
 $\Theta_r(\lambda)$ and $Z_r(\lambda)$ have the same properties. \end{Example}
 
 \bigskip
 
\section{Verma Modules of Complexity at most $2$} \label{s:appl}
The purpose of this final section is to indicate the utility of the results and methods established so far. In particular, we show that questions concerning Verma modules of complexity at most 
$2$ can be reduced to the case where the underlying groups are the first and second Frobenius kernels of $\SL(2)$. Throughout, $G$ is assumed to be defined over $\FF_p$.
 
 We begin with the following refinement of Lemma \ref{AR-V1}.
 
 \bigskip
 
 \begin{Lemma} \label{VMC1} Suppose that $G$ is semi-simple, simply connected and that $p$ is good for $G$. Let $Z_r(\lambda)$ be a Verma module of complexity $1$. Then 
 $\Theta_r(\lambda) \cong \ZZ[A_\infty]/(\tau)$ is a homogeneous tube. \end{Lemma}
 
\begin{proof} Writing $d+1 = \depth(\lambda)$, we obtain from Theorem \ref{DR6} the identity $r-d = 1$. Thus, Corollary \ref{DR7} provides an isomorphism
\[ \Theta_r(\lambda) \cong \Theta_1(\mu)\]
for some weight $\mu$. In particular, the component $\Theta_1(\mu)$ consists of $\tau$-periodic vertices, and \cite[(2.5)]{Fa1} yields the result. \end{proof} 
 
\bigskip
\noindent
As a first Corollary, we show that Verma modules of complexity $2$ possess equidimenisonal support varieties:
 
\bigskip
 
\begin{Corollary} \label{VMC2} Let $G$ be semi-simple and simply connected. Suppose that $p \ge 7$, and let $\lambda \in X(T)$ be a weight such that $\cx_{G_r}(Z_r(\lambda)) = 2$. 
Then $\cV_{G_r}(Z_r(\lambda))$ is equidimensional.\end{Corollary}
 
\begin{proof} In view of Theorem \ref{DR2}, we may assume that $\depth(\lambda) = 1$. Thanks to Corollary \ref{VV2}, we thus have
\[ 2 = \cx_{G_r}(Z_r(\lambda)) \ge r,\]
and Corollary \ref{DR5} yields the irreducibility of $\cV_{G_r}(Z_r(\lambda))$ in case $r=2$. Alternatively, \cite[(3.3)]{FR} implies that the rank variety $V_1(G)_{Z_1(\lambda)}$ is 
equidimenional. Owing to \cite[(1.4)]{FP1}, the support variety $\cV_{G_1}(Z_1(\lambda))$ enjoys the same property. \end{proof}

\bigskip
\noindent
By showing that Verma modules of complexity $2$ over higher Frobenius kernels essentially only occur for $\SL(2)_2$, our next result extends \cite[(3.5)]{FR}. We leave the modifications 
concerning modules of complexity $1$ to the interested reader.

\bigskip

\begin{Proposition} \label{VMC3} Suppose $r \ge 2$, $p \ge 7$, and let $\lambda \in X(T)$ be a weight such that $\depth(\lambda) = 1$ and $\cx_{G_r}(Z_r(\lambda)) \le 2$.  Then the 
following statements hold:

{\rm (1)} \ $r =2$ and $G_r \cong S_r \times H_r$, where $S \cong \SL(2)$ and $H$ is reductive.
 
{\rm (2)} \ There is an isomorphism $Z_r(\lambda) \cong Z_r^S(\lambda)\!\otimes_k\! Z_r^H(\lambda)$, with $\cx_{S_r}(Z_r^S(\lambda))= \cx_{G_r}(Z_r(\lambda))$, and the second 
factor being a simple projective $H_r$-module.  \end{Proposition} 

\begin{proof} Since $\Psi_\lambda \ne \Psi$ and $p$ is good, there exists a simple root $\alpha$ in $\Psi \setminus \Psi_\lambda$. According to Corollary \ref{VV2}, we have $r= 2$ as well 
as $\cV_{(U_\alpha)_r}(k) \subseteq \cV_{G_r}(Z_r(\lambda))$. In view of \cite[(1.14)]{SFB1} and  \cite[(5.2),(6.8)]{SFB2}, we obtain the same inclusion for the corresponding rank 
varieties: 
\[V_r(U_\alpha) \subseteq V_r(G)_{Z_r(\lambda)}.\]
Thus, $V_r(U_\alpha) \cong V_r(\GG_{a}) \cong k^r$ (cf.\ \cite[(1.10)]{SFB1}) is an irreducible component of the $B$-invariant variety $V_r(G)_{Z_r(\lambda)}$. Hence $V_r(U_\alpha)$ 
is $B$-invariant, and the orbit $B\dact \varphi$ of the canonical isomorphism $\GG_{a(r)} \stackrel{\varphi}{\lra} (U_\alpha)_r$ belongs to $V_r(U_\alpha)$. This implies that 
$(U_\alpha)_r$ is $B$-invariant, so that its Lie algebra $\Lie((U_\alpha)_r) = \Lie(U_\alpha) = \fg_\alpha$ enjoys the same property. An application of \cite[(3.3)]{FR} (and its proof) now 
implies that $G= SH$ is an almost direct product with 

(a) \ $e_\alpha \in \Lie (S)$, and

(b) \ $S$ is almost simple of  type $A_1$, and

(c) \ $\Lie(G) = \Lie(S)\oplus \Lie(H)$.

\noindent
Our assumption on $p$ implies $(S\cap H)_1 \subseteq {\rm Cent}(S)_1 = e_k$, so that $(S\cap H)_r = e_k$.  Thus, the canonical map $S\times H \lra G$ induces a closed embedding $S_r 
\times H_r \hookrightarrow G_r$. Applying \cite[(I.9.6)]{Ja3} several times we obtain
\[ \dim_k k[G_r] = p^{r\dim G} = p^{r(\dim S + \dim H)} = \dim_k k[S_r] \dim_k k[H_r] = \dim_k k[S_r\times H_r],\]
so that the above map is in fact an isomorphism. 

In view of \cite[(I.7.9(3))]{Ja3}, the arguments of \cite[(1.1)]{FR} show that the decomposition $G_r = S_r \times H_r$ induces an isomorphism 
\[ Z_r(\lambda) \cong Z_r^S(\lambda)\! \otimes_k \!Z_r^H(\lambda)\]
between $Z_r(\lambda)$ and the outer tensor product of the corresponding Verma modules of the factors. Thus, \cite[(7.2)]{SFB2} implies
\[ \cV_{G_r}(Z_r(\lambda)) \cong \cV_{S_r}(Z_r^S(\lambda))\! \times \!\cV_{H_r}(Z_r^H(\lambda)).\]
Since $e_\alpha \in \Lie (S)$, we see that $\alpha$ is a root of $S$. Consequently, $U_\alpha \subseteq S$, and another application of Corollary \ref{VV2} implies that 
$\cV_{(U_\alpha)_r}(k)  \subseteq \cV_{S_r}(Z_r^S(\lambda))$. In particular, $2 \le \dim \cV_{S_r}(Z_r^S(\lambda))$, so that $Z_r^H(\lambda)$ is projective. It thus follows from a 
consecutive application of \cite[(II.11.8)]{Ja3} and \cite[(II.9.6e)]{Ja3} that $Z_r^H(\lambda)$ is a simple projective $H_r$-module.

In view of \cite[(7.2.4)]{Sp}, the group $S$ is isomorphic to $\SL(2)$ or ${\rm PSL}(2)$. Since $p \ne 2$, the center of $\SL(2)$ is reduced, so that the canonical map $\SL(2) \lra {\rm 
PSL}(2)$ induces an isomorphism of the Frobenius kernels of these groups.\end{proof}

\bigskip
\noindent
Let $L \subseteq G$ be a Levi subgroup of $G$ and $S \subseteq  L$ be its semi-simple part. Given $\lambda \in X(T)$,  we denote by $\Theta^S_r(\lambda)$ the connected component of 
$\Gamma_s(S_r)$ containing the Verma module $Z_r^S(\lambda|_{T\cap S})$.

According to \cite[(II.9.1(3), II.9.4, II.11.8)]{Ja3} a non-projective Verma module $Z_r(\lambda)$ is defined by a highest weight $\lambda$ of depth $\depth(\lambda) \le r$. The following 
theorem reduces us to the case $r\le 2$. 

\bigskip

\begin{Theorem} \label{VMC4} Let $G$ be semi-simple and simply connected. Suppose that $p \ge 7$, and let $Z_r(\lambda)$ be a Verma module of complexity $\cx_{G_r}(Z_r(\lambda))= 
2$. Then we have $r':= r-\depth(\lambda)+1 \le 2$, and there exist a group $S \cong \SL(2), \SL(2)\!\times\!\SL(2), \SL(3)$, a weight $\gamma$ of  $S$ of depth $1$, and an isomorphism 
$\Theta_r(\lambda) \cong \Theta_{r'}(\gamma)$ sending $[Z_r(\lambda)]$ onto $[Z_{r'}(\gamma)]$. Moreover, we have $S \cong \SL(2)$ for $r'=2$.\end{Theorem}

\begin{proof} We put $\depth(\lambda) = d+1$, with $d \ge 0$. Thanks to Theorem \ref{DR6} and Corollary \ref{DR7}, there exists a weight $\mu$ of depth $1$ and a quiver isomorphism 
$\Theta_r(\lambda) \cong \Theta_{r-d}(\mu)$, sending $[Z_r(\lambda)]$ onto $[Z_{r-d}(\mu)]$. By Theorem \ref{DR6}, we have 
\[ r' := r-d \le \cx_{G_{r-d}}(Z_{r-d}(\mu)) = \cx_{G_r}(Z_r(\lambda)) = 2.\]
If $r' = 1$, then our result follows from \cite[(3.5)]{FR}. Alternatively, Proposition \ref{VMC3} implies the existence of a decomposition $G_{r'} = S_{r'}\! \times\! H_{r'}$ with $S \cong 
\SL(2)$. By the same token, there is an isomorphism 
\[ Z_{r'}(\mu) \cong Z_{r'}^S(\mu) \!\otimes_k\!Z_{r'}^H(\mu),\]
with $Z_{r'}^H(\mu)$ being a simple projective $H_{r'}$-module. By the arguments of \cite[(3.5)]{FR}, the functor 
\[ S_{r'}\text{-mod}  \lra G_{r'}\text{-mod} \ \ ; \ \ M \mapsto M\!\otimes_k\! Z_{r'}^H(\mu)\]
induces a Morita equivalence between the blocks $\cB_{r'}(\mu|_{T\cap S})$ and $\cB_{r'}(\mu)$ sending $Z_{r'}^S(\mu)$ onto $Z_{r'}(\mu)$. As a result, we obtain a quiver isomorphism 
$\Theta_{r'}^S(\mu) \cong \Theta_{r'}(\mu)$, which maps $[Z_{r'}^S(\mu)]$ onto $[Z_{r'}(\mu)]$.  Our assertion thus follows by composing the above isomorphisms of stable translation 
quivers. \end{proof}

\bigskip

\begin{Remark} By the same line of arguments, a semi-simple, simply connected group $G$ affords Verma modules of complexity $2$ only if it possesses almost simple factors of types $A_1$ 
or $A_2$. \end{Remark} 

\bigskip

\bigskip
\noindent
{\bf Acknowledgements}: This research was funded in part by EPSRC grant GR/R86546/01. Part of the research for this paper was carried out during a visit by the first author to the University of Southampton: he is grateful to the members of the School of Mathematics for their hospitality.

The authors would like to thank the referee for carefully reading the manuscript along with providing helpful suggestions.

\end{document}